\documentclass[11pt,reqno]{amsart}
\usepackage{amssymb, amsmath, amsthm,amsrefs,marvosym,wasysym}
\usepackage{latexsym}
\usepackage{color}
\usepackage{exscale}
\usepackage{fullpage}
\usepackage{rotating}
\usepackage{bm, bbm, mathrsfs}
\usepackage{graphics, graphicx}
\newtheorem{theorem}{Theorem}
\newtheorem{lemma}[theorem]{Lemma}
\newtheorem{proposition}[theorem]{Proposition}

\newtheorem{definition}{Definition}[section]

\newtheorem{remark}{Remark}[section]

\def\eq#1{(\ref{#1})}

\def\({\left(\begin{array}{cccccc}}
\def\){\end{array}\right)}

\def\eq#1{(\ref{#1})}

\def\({\left(\begin{array}{cccccc}}
\def\){\end{array}\right)}

\def\bes{\begin{eqnarray}}
\def\ees{\end{eqnarray}}

\newcommand{\beq}{\begin{equation}}
\newcommand{\eeq}{\end{equation}}
\newcommand{\bea}{\begin{eqnarray}}
\newcommand{\eea}{\end{eqnarray}}
\newcommand{\beann}{\begin{eqnarray*}}
\newcommand{\eeann}{\end{eqnarray*}}

\newcommand{\eps}{\ensuremath{\varepsilon}}
\newcommand{\ep}{\ensuremath{\epsilon}}

\newcommand{\lam}{\ensuremath{\lambda}}

\newcommand{\RR}{\mathbb{R}}

\DeclareMathOperator{\sgn}{sgn}

\newcommand{\lj}{\big[\!\!\big[}
\newcommand{\rj}{\big]\!\!\big]}

\numberwithin{equation}{section}

\begin{document}

\title{1-D Isentropic Euler flows: Self-similar Vacuum Solutions}

\begin{abstract}
	We consider one-dimensional self-similar solutions to
	the isentropic Euler system when the initial data 
	are at vacuum to the left of the origin. For $x>0$ the initial velocity 
	and sound speed are of form 
	$u_0(x)=u_+x^{1-\lambda}$ and $c_0(x)=c_+x^{1-\lambda}$,
	for constants $u_+\in\RR$, $c_+>0$, $\lambda\in\RR$.
	We analyze the resulting solutions in terms of 
	the similarity parameter $\lambda$, the adiabatic exponent 
	$\gamma$, and the initial (signed) Mach number $\text{Ma}=u_+/c_+$. 
	
	Restricting attention to locally bounded data, we find that 
	when the sound speed initially decays to zero in a H\"older manner
	($0<\lambda<1$), the resulting flow is always defined globally.
	Furthermore, there are three regimes depending on $\text{Ma}$:
	\begin{itemize}
		\item for sufficiently large positive $\text{Ma}$-values, the solution is continuous
		and the initial H\"older decay is immediately replaced by $C^1$-decay to 
		vacuum along a stationary vacuum interface;
		\item for moderate values of $\text{Ma}$, the solution is again continuous and with 
		an accelerating vacuum interface along which $c^2$ decays linearly to zero
		(i.e., a ``physical singularity'');
		\item for sufficiently large negative $\text{Ma}$-values, the solution contains 
		a shock wave emanating from the initial vacuum interface and propagating 
		into the fluid, together with a physical singularity along an accelerating vacuum 
		interface. 
	\end{itemize}

	In contrast, when the sound speed initially decays to zero in a $C^1$ manner 
	($\lambda<0$), a global flow exists only for sufficiently large positive values of $\text{Ma}$. 
	Non-existence of global solutions for smaller $\text{Ma}$-values is due to rapid growth 
	of the data at infinity and is unrelated to the presence of a vacuum.
\end{abstract}

\author{Helge Kristian Jenssen}
\thanks{Address: Department of
Mathematics, Penn State University,
State College, PA 16802, USA\\
Email: jenssen@math.psu.edu\\
ORCHID: 0000-0002-5344-4399\\
Acknowledgements:
This material is based in part upon work supported by the National Science Foundation 
under Grant Number DMS-1813283. Any opinions, findings, and conclusions 
or recommendations expressed in this material are those of the authors and do not 
necessarily reflect the views of the National Science Foundation}

\date{\today}
\maketitle

%

\tableofcontents

\section{Introduction}\label{intro}
The compressible Euler system \eq{mass}-\eq{mom} degenerates 
at vacuum (losing strict hyperbolicity), and already local 
existence of solutions with vacuum is a non-trivial issue. 
The analysis of how vacuum interfaces propagate in isentropic flow, 
including a precise description of decay to vacuum, has been addressed 
in a number of works in recent years. 
One part of this effort addresses the propagation of a so-called ``physical singularity''
in the sound speed along the vacuum interface. The seminal work \cite{liu} argued that,
generically, the sound speed suffers a square root singularity whenever the 
interface is accelerated by the internal pressure. A series of recent works has 
demonstrated local existence and stability of solutions with a physical singularity,
both in one and several space dimensions. 

A related line of inquiry concerns the situation when the initial decay of the sound 
speed is different from that of the physical singularity. For concreteness, consider a
1-dimensional situation with the fluid initially located to the right of the origin, and suppose
the initial sound speed decays to zero like $x^{1-\lambda}$ as $x\downarrow0$. 
It has been conjectured that any H\"older decay ($0<\lambda<1$) should immediately
lead to acceleration of the interface, together with an instantaneous switch to a physical
singularity. Furthermore, smooth decay to vacuum ($\lambda<0$)
should generically be replaced by a physical singularity after a finite waiting time, during 
which the pressure builds up. 

In this work we study the particular class of self-similar isentropic flows 
with vacuum, obtaining concrete examples of how vacuum interfaces propagate.
Beside their intrinsic interest they provide relevant insights about the 
above conjectures. 

First, these particular solutions demonstrate that, with a fixed 
decay of the sound speed to vacuum, different initial velocities can yield distinct 
qualitative behaviors (a stationary interface along which $c$ decays 
smoothly to vacuum vs.\ an accelerating interface with a physical singularity).

A second, more striking, feature is the possibility
of a shock wave emanating from the initial interface and moving into the fluid. 
It turns out that, for the special solutions under consideration, once 
the adiabatic constant $\gamma$ and the similarity parameter $\lambda$ are fixed,
the only relevant parameter is the (signed) Mach number $\text{Ma}$
of the initial data. 
For self-similar solutions with a vacuum on the left, we show that a shock is 
necessarily generated whenever the initial Mach number is sufficiently large 
and negative, i.e., whenever the gas initially moves sufficiently fast toward the vacuum. 

After describing the setup for self-similar vacuum flow in Section \ref{ss_euler},
we state our main findings in Section \ref{results}.
For conciseness we provide these in the case with $0<\lambda<1$ and $1<\gamma<3$; 
the (similar) conclusions for other cases are stated later.
Section \ref{interpretation} provides a physical interpretation of our findings,
followed by a discussion of related works in Section \ref{other_works}.
Section \ref{sim_odes} records the similarity ODEs together with an 
outline of the construction. 
Sections \ref{P0_P2} and \ref{P3_P6} provide the relevant details about the 
(up to seven) critical points of the similarity ODEs. Several of these 
are located along two special, straight-line trajectories $E_\pm$ that 
are used later to delimit distinct types of solution behaviors.
Section \ref{ss_grps_rh_e_conds} contains an analysis of the 
Rankine-Hugoniot relations and entropy conditions for self-similar flows, 
including properties of the ``Hugoniot locus'' corresponding to a given trajectory
of the similarity ODEs.
These results are then used in Sections \ref{I}-\ref{II} to build physically 
admissible self-similar vacuum flows when the adiabatic constant satisfies
$1<\gamma<3$. The corresponding results for $\gamma\geq 3$ are 
qualitatively the same, and these are described in Section \ref{other_cases}.

\section{Self-similar Euler flows}\label{ss_euler}
The 1-d isentropic compressible Euler system expresses conservation of 
mass and linear momentum in isentropic flow of an ideal gas with planar symmetry: 
\begin{align}
	\rho_t+(\rho u)_x&=0 \label{mass}\\
	(\rho u)_t+(\rho u^2)_x+p_x&=0.\label{mom}
\end{align}
Here the independent variables are time $t$ and position $x\in\RR$, and the 
primary dependent variables are density $\rho(t,x)$ and fluid velocity $u(t,x)$.
The pressure $p$ given by
\beq\label{pressure1}
	p(\rho)=a^2\rho^\gamma,
\eeq
where the adiabatic constant satisfies $\gamma>1$, and $a>0$ is a constant.
The local speed of sound $c\geq 0$ is given by
\beq\label{sound_speed}
	c=\sqrt{p'(\rho)}=a\sqrt{\gamma}\rho^\frac{\gamma-1}{2}.
\eeq
In terms of $u$ and $c$ the system takes the following form in smooth regions of the flow:
\begin{align}
	u_t+ uu_x +\ell cc_x&= 0\label{u}\\
	c_t+uc_x+\ell^{-1} cu_x&=0,\label{c}
\end{align}
where
\[\ell:=\textstyle\frac{2}{\gamma-1}.\]
The initial data are 
\[u(0,x)=u_0(x)\qquad\text{and}\qquad c(0,x)=c_0(x)\geq0.\]
We are interested in the particular class of initial data that generate 
self-similar solutions to \eq{u}-\eq{c}.

\subsection{Self-similar flows}\label{ss_flows}
For a given solution $(u,c)$ of \eq{u}-\eq{c}, and any $\ep>0$ and $\lambda\in\RR$,
\[u_\ep(t,x):=\eps^{1-\lambda} u(\textstyle\frac{t}{\ep^\lambda},\textstyle\frac{x}{\ep}),\qquad
c_\ep(t,x):=\eps^{1-\lambda} c(\textstyle\frac{t}{\ep^\lambda},\textstyle\frac{x}{\ep})\]
is again a solution of \eq{u}-\eq{c}. The given solution is {\em self-similar}
provided $(u,c)\equiv(u_\ep,c_\ep)$ for all $\ep>0$, i.e.,
\beq\label{ss_soln}
	u(t,x)=\eps^{1-\lambda} u(\textstyle\frac{t}{\ep^\lambda},\textstyle\frac{x}{\ep}), \qquad
	c(t,x)=\eps^{1-\lambda} c(\textstyle\frac{t}{\ep^\lambda},\textstyle\frac{x}{\ep})
	\qquad\text{for all $t,x$ and all $\ep>0$.}
\eeq
Evaluating these at $(t,x)=(0,\pm1)$ and setting $y=\ep^{-1}$,
we get that the initial data $u_0$, $c_0$ satisfy
\[u_0(\pm y)=y^{1-\lambda} u_0(\pm1),\qquad c_0(\pm y)=y^{1-\lambda} c_0(\pm1)\qquad
\text{for any $y>0$.}\]
This shows that the initial data for a self-similar solution of \eq{u}-\eq{c} must 
be proportional to $|x|^{1-\lambda}$, possibly with different constants of proportionality 
for $x\gtrless 0$. 

Evaluating \eq{ss_soln} with $\ep=t^\frac{1}{\lambda}$ shows that 
a self-similar solution of \eq{u}-\eq{c} has the form
\[u(t,x),c(t,x)=t^{\frac{1}{\lambda}-1}\times \text{[Function of $t^{-\frac{1}{\lambda}}x$].}\]
Following \cites{cf,laz} we opt to use the similarity coordinate 
\beq\label{sim_coord}
	\xi:=t^{-\frac{1}{\lambda}}x,
\eeq
and posit 
\beq\label{sim_var_u}
	u(t,x)=-\textstyle\frac{1}{\lambda}\frac{x}{t}V(\xi),
\eeq
\beq\label{sim_var_c}
	c(t,x)=-\textstyle\frac{1}{\lambda}\frac{x}{t}C(\xi).
\eeq
Alternatively, as $t>0$ so that $\sgn(\xi)=\sgn(x)$, we have
\begin{align}
	u(t,x)
	&=-\textstyle\frac{1}{\lambda}t^{\frac{1}{\lambda}-1}\xi V(\xi) 
	=-\textstyle\frac{1}{\lambda}\sgn(\xi)|\xi|^\lambda V(\xi)|x|^{1-\lambda},
	\label{sim_var_u_alt}\\
	c(t,x)
	&=-\textstyle\frac{1}{\lambda}t^{\frac{1}{\lambda}-1}\xi C(\xi)
	=-\textstyle\frac{1}{\lambda}\sgn(\xi)|\xi|^\lambda C(\xi)|x|^{1-\lambda}.
	\label{sim_var_c_alt}
\end{align}
The variables $V$ and $C$
satisfy a coupled system of similarity ODEs recorded in 
\eq{V_ode}-\eq{C_ode} below.

\begin{remark}\label{c_sign}
	Both $V(\xi)$ and $C(\xi)$ may be of either sign. However,
	as we consider flows defined for positive times,
	the convention $c(t,x)\geq0$ imposes
	(according to \eq{sim_var_c_alt}) the constraint  
	\beq\label{sign_cond}
		\textstyle\frac{\xi}{\lambda} C(\xi)\leq 0 \qquad\text{for all $\xi\in\RR$.}
	\eeq
\end{remark}

The analysis above shows that it is natural to consider initial value problems 
for the 1-d isentropic Euler system \eq{u}-\eq{c} 
with data of the form $u_0(x)=u_\pm|x|^{1-\lambda}$ and 
$c_0(x)=c_\pm|x|^{1-\lambda}$ for $x\gtrless0$,
where $u_\pm\in\RR$ and $c_\pm\geq 0$ are constants.

The case $\lambda=1$ corresponds to standard Riemann problems where 
the data consist of two constant states separated by a jump discontinuity.
The solution in this case is well known \cites{sm,gr}, including the case 
where one of the constant states is at vacuum \cites{gb,st}.
The case $\lambda=0$ yields piecewise linear initial data for $u$ and $c$; the 
formulation above is not appropriate in this case, and we leave it out of the following 
discussion. 
Also, $\lambda>1$ yields unbounded data; while such data can still
have locally bounded energy (specifically, for $\lambda<\frac{3\gamma-1}{2\gamma}$), 
we find it more relevant to consider data with locally bounded amplitudes. 
We thus restrict to
\begin{itemize}
	\item{} Case (I): $0<\lambda<1$; $u_0$, $c_0$ are locally bounded 
	and $C^{0,1-\lambda}$; or
	\item{} Case (II): $\lambda< 0$; $u_0$, $c_0$ are locally bounded 
	and $C^1$.
\end{itemize}
Our objective is to exploit the particular class of self-similar Euler flows to
study solutions with a vacuum region, chosen to be located initially to the 
left of $x=0$. To this 
end we further restrict attention to data for which the initial density 
(equivalently, sound speed) vanishes for $x<0$:
\beq\label{vac_ss_data}
	u_0(x)=\left\{\begin{array}{ll}
	\text{undefined} & x<0\\
	u_+x^{1-\lambda} & x>0,
	\end{array}\right.
	\qquad
	c_0(x)=\left\{\begin{array}{ll}
	0 & x<0\\
	c_+x^{1-\lambda} & x>0,
	\end{array}\right.
\eeq 
where $u_+\in\RR$ and $c_+> 0$ are constants. The goals are to describe 
the behavior of the resulting solution in terms of how the vacuum interface 
moves, how $c$ decays to zero along it at fixed times, and whether shock 
waves are present or not. 
In particular, we are interested in identifying any instantaneous change in 
the decay rate of $c$ from $t=0$ to $t=0+$.
This requires a detailed study of the phase portrait of the similarity ODEs 
\eq{V_ode}-\eq{C_ode}, as well as the jump and admissibility conditions along
discontinuities.

\section{Main results}\label{results}
We use the following standard terminology. At a fixed time and with leading order 
behavior of the sound speed near the vacuum interface given by 
$c\sim [\text{distance to interface}]^\alpha$, we say that $c$ decays 
to vacuum in a $C^1$ (H\"older, Lipschitz) manner provided $\alpha>1$ 
($0<\alpha<1$, $\alpha=1$, respectively). A ``physical singularity'' 
refers to the case $\alpha=\frac{1}{2}$. 
Also, the signed Mach number of the data \eq{vac_ss_data} is denoted 
\beq\label{mach}
	\text{Ma}:=\textstyle\frac{u_+}{c_+}.
\eeq

\subsection{Case (I): $0<\lambda<1$ and $1<\gamma<3$}
We next formulate our main findings for the case where the initial sound 
speed decays to vacuum in a H\"older manner ($0<\lambda<1$),
and for the physically more relevant case $1<\gamma<3$.
It turns out that all the key features are present in this case, with only 
minor differences when $\gamma\geq 3$ (see Section \ref{other_cases}). On the other hand, for $C^1$ 
decay to vacuum ($\lambda<0$) we obtain a global self-similar solution 
only for a restricted range of initial $\text{Ma}$-numbers; see Section 
\ref{interpretation} for a discussion of this issue. (Recall that $\ell=\frac{2}{\gamma-1}$.) 

\begin{theorem}\label{thm_1}
	Assume $1<\gamma<3$ and $0<\lambda<1$. Then the initial 
	value problem \eq{u}-\eq{c}-\eq{vac_ss_data} admits a self-similar 
	solution $(u(t,x),c(t,x))$ of the form 
	\eq{sim_coord}-\eq{sim_var_u}-\eq{sim_var_c} defined for all $t>0$, $x\in\RR$. 
	Furthermore, we have the following behaviors depending on the initial Mach-number:
		\begin{itemize} 
			\item[(a)] For $\mathrm{Ma}> \ell$ the solution is continuous and 
			with a stationary vacuum interface at $\{x=0\}$, along which the sound 
			speed $c$ decays to vacuum in a $C^1$ manner at each time $t>0$. 
			The same applies to the limiting case $\text{Ma}=\ell$, except that 
			$c$ then decays to vacuum in a Lipschitz manner at each time $t>0$. 
			\item[(b)] For $-\ell\leq\mathrm{Ma}<\ell$ the solution is continuous 
			and with a non-stationary vacuum interface propagating to the left along 
			$x=t^\frac{1}{\lambda}\xi_v$ for some $\xi_v<0$. Along the interface the 
			sound speed $c$ exhibits a physical singularity at each time $t>0$.
			\item[(c)] For $\mathrm{Ma}<-\ell$ the solution is discontinuous and 
			contains a single 2-shock, together with a non-stationary vacuum 
			interface. The vacuum interface propagates to the left along 
			$x=t^\frac{1}{\lambda}\xi_v$ for some $\xi_v<0$, while the 2-shock 
			propagates to the left along $x=t^\frac{1}{\lambda}\xi_s$, where 
			$\xi_v<\xi_s<0$. Along the interface the sound speed $c$ exhibits a 
			physical singularity at each time $t>0$.
		\end{itemize}
		In particular, for all cases except (b) and (c) with $\lambda=\frac{1}{2}$, 
		there is an abrupt change from $t=0$ to $t=0+$ in the decay of 
		$c(t,\cdot)$ to vacuum.
\end{theorem}

\begin{remark}
	While we expect that any globally defined and self-similar solution 
	of the initial value problem \eq{u}-\eq{c}-\eq{vac_ss_data} is unique, 
	we have not been able to prove this in all cases under consideration.
	Specifically, the occurrence of a shock wave is established 
	by demonstrating that a certain trajectory of the similarity ODEs 
	intersects the Hugoniot locus of a certain other trajectory. While an 
	argument based on continuity shows that an intersection must occur 
	in the relevant cases, we do not have a proof that the intersection 
	is unique.
\end{remark}

\subsection{Physical interpretation}\label{interpretation}
With the initial vacuum region along $\{x<0\}$, $u_+>0$ ($u_+<0$) 
corresponds to the situation where 
the fluid initially moves away from (toward) the vacuum.
To interpret physically our findings
we note the following about the initial data \eq{vac_ss_data}:
\begin{itemize}
	\item[-] the initial density and pressure increase as we move 
	to the right from the interface at $x=0$;
	\item[-] the same applies to the initial speed, i.e., $|u_0(x)|$;
	\item[-] the fluid is initially at rest at $x=0$.
\end{itemize}
The initial pressure gradient tends to push the fluid to the left. 
With the fluid initially at rest at $x=0$, this situation is conducive
to shock formation
- a fortiori if the fluid is initially moving toward the vacuum 
($u_+<0$). On the other hand, if $u_+>0$ then the 
initial velocity $u_0(x)$ increases with $x$. 
Clearly, this tends to rarify the fluid and acts
against shock formation. Finally, the fluid is free to expand into the vacuum 
on the left, and this also provides a rarefying effect. In particular, one may 
expect that this last effect can, to some degree, counteract the tendency 
of shock formation due to an initial leftward velocity field.

For data \eq{vac_ss_data} with $c_+$ fixed, it is therefore reasonable 
to expect shock formation when the fluid is initially moving sufficiently fast
toward the vacuum. Our results show that this is indeed the case and also 
provide a concise criterion for what ``sufficiently fast'' means. 

Consider case (I) ($0<\lambda<1$, $c$ decays to vacuum in a H\"older manner), in which the initial fluid speed $|u_0(x)|$
and sound speed $c_0(x)$ grow sub-linearly as $x$ increases. 
In this case we find that the data \eq{vac_ss_data} always give a globally 
defined Euler flow. Furthermore, the solution is shock-free if and only if the initial 
Mach number satisfies
\beq\label{soln_restricn}
	\text{Ma}\geq \text{Ma}_{cr},
\eeq
where the critical Mach number $\text{Ma}_{cr}<0$ in general depends on 
both $\lambda$ and $\gamma$. (For the case $1<\gamma<3$ treated in 
Theorem \ref{thm_1}, we have $\text{Ma}_{cr}=-\ell$ independently of $\lambda$;
for $\gamma>3$, $\text{Ma}_{cr}$ also depends on $\lambda$.) 
When \eq{soln_restricn} is violated the initial pressure gradient ``wins,'' and a single
shock is necessarily generated at initial time $t=0+$, emanating
from the initial vacuum interface at $x=0$. It propagates into the fluid (either 
to the left or to the right, depending on $\gamma$, $\lambda$, $\text{Ma}$), 
i.e., its distance to the vacuum interface increases with time. Furthermore, 
in our setup with the vacuum region located 
to the left of the fluid region, any shock is necessarily a 2-shock.
Thus, as time increases, fluid particles enter the region between the shock 
and the vacuum interface.
\begin{remark}\label{vac_reg}
        It is interesting to compare this 
        behavior with the situation where a shock  approaches a vacuum interface
        from within the fluid, see \cites{sak,ls}. Such a shock weakens
        and disappears as it reaches the vacuum. In this sense, a vacuum provides 
        a regularizing effect. Our results add to the picture by showing that 
        {\em a shock wave can be generated at a vacuum interface}. 
\end{remark}

For case (II) ($\lambda<0$, $c$ decays to vacuum in a $C^1$ manner) 
the situation is different and less complete (see Section \ref{II}). 
Again, we shall find that data with a sufficiently large $\text{Ma}$-value
generate a globally defined and shock-free solution. 
However, for smaller $\text{Ma}$-values we no longer 
obtain a globally defined flow. Instead, we will be able to continue the solution 
in a continuous manner only to a region in the $(x,t)$-plane of the form 
$\{(x,t)\,|\, 0<t^{1/|\lambda|}x<\xi_*\}$.
Here, $\xi_*$ is a positive and finite number beyond which the underlying 
solution of the similarity ODEs cannot be continued because it runs into a 
singularity located along the critical (sonic) line $L_+$, cf.\ \eq{crit_lines}. 

This failure to resolve all initial value problems in case (II) is unsurprising
and not related to the presence of the vacuum. Indeed, in case (II) the initial sound 
speed $c_0(x)$ grows super-linearly, implying rapid growth of the pressure 
as $x\uparrow\infty$. Unless the initial velocity is sufficiently large and positive
it is reasonable to expect that an 
unbounded wave is immediately generated at $t=0+$, moving in from 
$x=+\infty$ along a path with $t^{-1/\lambda}x=\xi_*$ and leaving the flow 
undefined in its wake.
For this reason we do not analyze case (II) further when \eq{soln_restricn} is violated.

\subsection{Connections to other works}\label{other_works}
The analysis of how vacuum interfaces propagate in Euler flow, 
including a precise description of decay to vacuum, has been addressed 
in a number of works in recent years. 
An important motivation for the study of vacuum interfaces appears 
in connection with gaseous stars where pressure and gravity balance
\cites{ma,ni}. The early works \cites{st,gb} 
describe how a gas initially at rest and at constant pressure expands 
into a vacuum: Along the interface, which in this case propagates 
at constant speed, the sound speed $c$ decays to vacuum in a Lipschitz 
manner and connects continuously to the unperturbed state upstream.

In \cite{liu} Liu argued that when a vacuum interface propagates with 
non-zero acceleration, it is instead the quantity $c^2$ that decays
linearly to vacuum. Thus, $c$ should generically suffer a square-root 
singularity along a vacuum interface, a behavior referred to as a ``physical singularity.''
(As far as we know, this was first identified as the relevant boundary condition 
along an accelerating vacuum interface by Richtmyer and Lazarus \cite{rl}.)
By now, local existence of solutions propagating an initial physical singularity
has been established in both 1-d and multi-d 
\cites{jm1,jm2,cs1,cs2,dit_22,it_20,jm}.

Next consider a general vacuum initial value problem for \eq{u}-\eq{c}: 
\beq\label{gen_vac_ss_data}
	u_0(x)=\left\{\begin{array}{ll}
	\text{undefined} & x<0\\
	\bar u(x) & x>0,
	\end{array}\right.
	\qquad
	c_0(x)=\left\{\begin{array}{ll}
	0 & x<0\\
	\bar c(x) & x>0.
	\end{array}\right.
\eeq 
In order to compare with the particular case analyzed in the present work
(i.e., $\bar u(x)=u_+x^{1-\lambda}$, $\bar c(x)=c_+x^{1-\lambda}$, cf.\
\eq{vac_ss_data}), we assume that the initial sound speed in \eq{gen_vac_ss_data}
satisfies $\bar c(x)\sim x^{1-\lambda}$ as $x\downarrow 0$.
Here, $\lambda=\frac{1}{2}$ corresponds to the presence of an initial physical singularity.
For $\lambda\neq\frac{1}{2}$, the following conjectures have been formulated
\cites{ly,jm_12}:
\begin{itemize}
	\item[(a)] For $0<\lambda<1$, $\lambda\neq\frac{1}{2}$, the interface 
	immediately accelerates and propagates a physical singularity.
	\item[(b)] For $\lambda<0$ local existence 
	of a solution propagating 
	the same decay, is known \cites{ly1,jm}; it is conjectured that after a finite 
	{\em waiting time}, during which the pressure builds up, the interface 
	starts to accelerate and a physical singularity appears.
\end{itemize}
Note that (a) and (b) correspond to our cases (I) and (II), respectively.
However, as the present work deals with non-generic and unbounded data, 
our setup does not fully correspond to the setting of these conjectures. 
E.g., as noted above, for case (II) we obtain a physically meaningful solution 
only when the gas is initially moving sufficiently fast away from the vacuum.
When this is the case, the resulting interface is stationary and the ``waiting time'' may be 
said to be infinite. There is no essential contradiction with (b) in this.
We note that a special (and, as far as we know, first) 
instance of waiting time behavior in Euler flow is exhibited in the recent work
\cite{cfopp} on dam breaking.

Also our results for case (I) are essentially in agreement with (a), again 
with the caveat that the interface may remain stationary 
due to a sufficiently strong rarefying initial velocity field (in which case 
no physical singularity appears). In all other cases
we have immediate acceleration and propagation of a physical singularity
along the interface, in accordance with (a).

What the examples of self-similar Euler flows do demonstrate in 
connection with the conjectures above, are two things. First, the  
velocity field is as important as the sound speed in determining the 
qualitative features of the resulting flow. And second, when $c$ and $u$
decay to vacuum in a H\"older manner, a single 
shock is immediately generated at the interface whenever the gas 
initially moves sufficiently fast toward the vacuum. This provides a concrete 
case of the ``ill-posedness'' described in \cite{jm_12}. In particular, any
existence theory for isentropic Euler flows which covers data with 
general decay to vacuum, must cover the possibility that a shock wave 
is generated at the initial interface.

\section{Similarity ODEs and outline}\label{sim_odes}
We start by recording the similarity ODEs that describe self-similar Euler flows,
together with some of their basic properties. These are then used to 
provide an outline of the construction of the trajectories needed to resolve
the initial value problem \eq{u}, \eq{c}, \eq{vac_ss_data}.

\subsection{Similarity ODEs}
Substitution of \eq{sim_var_u}-\eq{sim_var_c} into \eq{u}-\eq{c} yields 
\begin{align}
	\frac{dV}{d\xi}&=\frac{1}{\xi}\frac{G(V,C)}{D(V,C)}\label{V_ode}\\
	\frac{dC}{d\xi}&=\frac{1}{\xi}\frac{F(V,C)}{D(V,C)},\label{C_ode}
\end{align}
where 
\begin{align}
	D(V,C)&=(1+V)^2-C^2\label{D}\\
	G(V,C)&=C^2(V-V_*)-V(1+V)(\lam+V)\label{G}\\
	F(V,C)&=C\left\{C^2-(1+V)^2+k_1(1+V)-k_0\right\},\label{F}
\end{align}
with
\beq\label{V_*}
	V_*=\ell(1-\lambda)
\eeq
and
\beq\label{ks}
	k_1=(\ell^{-1}-1)(\lambda-1),\qquad
	k_0=\ell^{-1}(\lam-1).
\eeq
We record the symmetries
\beq\label{symms}
	G(V,-C)=G(V,C),\qquad F(V,-C)=-F(V,C),
\eeq
and also the fact that 
\beq\label{non_obvious_reln}
	F(V,\pm(1+V))\equiv \mp\ell^{-1}G(V,\pm(1+V)).
\eeq

The trajectories of the similarity ODEs provide values of $V(\xi)$ and $C(\xi)$, 
which in turn provide values of the flow variables $u(t,x)$, $c(t,x)$, 
via \eq{sim_var_u}-\eq{sim_var_c}, along the curves $\xi=t^{-1/\lambda}x\equiv 
constant$ in the $(x,t)$-plane.
The latter curves foliate the $(x,t)$-plane, and thus provide (at least in principle) a 
complete Euler flow defined for all $t>0$.

The key to analyzing self-similar Euler flows is the fact that 
\eq{V_ode}-\eq{C_ode} yields the single {\em reduced similarity ODE}
\beq\label{CV_ode}
	\frac{dC}{dV}=\frac{F(V,C)}{G(V,C)}
\eeq
relating $V$ and $C$ along similarity solutions. Trajectories 
(solution curves) of the original similarity ODEs 
\eq{V_ode}-\eq{C_ode} can therefore be analyzed  via the 
phase portrait of \eq{CV_ode}.

Some care is required: while closely related, the ODE system \eq{V_ode}-\eq{C_ode} 
and the single ODE \eq{CV_ode} are not equivalent.
Specifically, to obtain physically meaningful solutions to the Euler system 
we need solutions to the original ODE system \eq{V_ode}-\eq{C_ode}.
However, in contrast to the reduced ODE \eq{CV_ode}, the system 
is singular along the {\em critical (sonic) lines}
\beq\label{crit_lines}
	L_\pm:=\{C=\pm(1+V)\},
\eeq
across which $D(V,C)$ changes sign. According to \eq{non_obvious_reln},
$F(V,C)$ and $G(V,C)$ are proportional along $L_\pm$.
Therefore, if a trajectory $\Gamma$ of \eq{CV_ode} crosses one of the critical lines 
$L_\pm$ at a point $P$ where one, and hence both, of $F$ and $G$ 
are non-zero, then the flow of \eq{V_ode}-\eq{C_ode} along $\Gamma$ is directed in 
opposite directions on either side of the critical line at $P$. This renders the trajectory 
$\Gamma$ useless for constructing a physically meaningful Euler flow from it. 

The upshot is a reduction of the set of relevant trajectories: any continuous crossing 
of one of the critical lines $L_\pm$ must occur at a ``triple point'' where $F$, $G$, 
$D$ all vanish. As we shall see in Section \ref{P3_P6}, there are 
always at least two triple points present (except in the particular case $\gamma=3$). 

We stress that this reduction is less drastic than what might first appear. This is because 
the critical lines can (indeed, must) be crossed by jumping from one trajectory
to another whenever the jump corresponds to an entropy admissible shock in 
the Euler flow. This issue is analyzed in Section \ref{ss_grps_rh_e_conds}.

\subsection{Outline of construction}\label{outline}
Taking uniqueness of self-similar solutions to \eq{u}-\eq{c}-\eq{vac_ss_data}
for granted, the solution is characterized by the requirement that, 
for each fixed $x>0$,
\beq\label{init_rates}
	\lim_{t\downarrow0} u(t,x)=u_+ x^{1-\lambda},\qquad
	\lim_{t\downarrow0} c(t,x)=c_+x^{1-\lambda}.
\eeq
Also, the choice $\xi:=t^{-\frac{1}{\lambda}}x$ for the similarity variable 
implies that
\[t\downarrow 0\quad\text{with $x>0$ fixed, corresponds to}\quad
\left\{\begin{array}{ll}
	\xi\downarrow 0 & \text{when $\lambda<0$,}\\ 
	\xi\uparrow\infty & \text{when $\lambda>0$.}
\end{array}\right.\]
Therefore, with the initial data \eq{vac_ss_data}, the relations 
\eq{sim_var_u_alt}-\eq{sim_var_c_alt} show that the corresponding 
solution $(V(\xi),C(\xi))$ of the similarity ODEs must satisfy
\beq\label{t=0_1}
	\lim_{\xi\downarrow0} \xi^\lambda V(\xi)=-\lambda u_+\qquad\text{and}\qquad  
	\lim_{\xi\downarrow0} \xi^\lambda C(\xi)=-\lambda c_+\qquad\text{when $\lambda<0$,}
\eeq
and 
\beq\label{t=0_2}
	\lim_{\xi\uparrow\infty} \xi^\lambda V(\xi)=-\lambda u_+\qquad\text{and}\qquad  
	\lim_{\xi\uparrow\infty} \xi^\lambda C(\xi)=-\lambda c_+\qquad\text{when $\lambda>0$.}
\eeq
In particular, we have
\beq\label{VC_at_t=0}
	(V(\xi),C(\xi))\to(0,0) \qquad
	\left\{\begin{array}{l}
		\text{as $\xi\downarrow0$ when $\lambda<0$, or }\\
		\text{as $\xi\uparrow\infty$ when $\lambda>0$.}
	\end{array}\right.
\eeq
Next, we record the fact (see Section \ref{P0}) that the origin $P_0=(0,0)$ 
in the $(V,C)$-plane is a critical point for the reduced similarity ODE \eq{CV_ode}. 
Furthermore, for any values of the parameters $\gamma$ and $\lambda$, $P_0$ 
is a star point (proper node), i.e., the linearization of \eq{CV_ode} there is 
$dC/dV=C/V$. A standard ODE result (\cite{ha}, Theorem 3.5 (iv)) implies 
a one-to-one correspondence between trajectories of \eq{CV_ode} 
approaching the origin, and their limiting 
slopes there. It follows from \eq{t=0_1}, \eq{t=0_2}, and \eq{VC_at_t=0}, that this limiting slope is 
given by the prescribed data \eq{vac_ss_data} according to
\beq\label{slope_at_t=0}
	\frac{C(\xi)}{V(\xi)} \to \frac{c_+}{u_+}\equiv \frac{1}{\text{Ma}}
	\qquad
	\left\{\begin{array}{l}
		\text{as $\xi\downarrow0$ when $\lambda<0$, or }\\
		\text{as $\xi\uparrow\infty$ when $\lambda>0$.}
	\end{array}\right.
\eeq
Finally, by combining \eq{VC_at_t=0} with the constraint in \eq{sign_cond} 
we obtain that the trajectory $(V(\xi),C(\xi))$ approaches the origin:
\beq\label{approach_at_t=0}
	\left\{\begin{array}{l}
		\text{with $C(\xi)>0$ as $\xi\downarrow0$ when $\lambda<0$, and}\\
		\text{with $C(\xi)<0$ as $\xi\uparrow\infty$ when $\lambda>0$.}
	\end{array}\right.
\eeq
The strategy for building the flow with the initial data \eq{vac_ss_data}
can now be made more precise: 
\begin{itemize}
	\item The data \eq{vac_ss_data} selects, via \eq{slope_at_t=0}, the relevant 
	trajectory of the similarity ODEs \eq{V_ode}-\eq{C_ode} near the origin in the 
	$(V,C)$-plane. This trajectory is denoted $\Gamma_0$; it lies in the upper 
	or lower half-plane for $\lambda<0$ or $\lambda>0$, respectively.
	\item As the solution of \eq{V_ode}-\eq{C_ode} moves away from the origin 
	along $\Gamma_0$ (with $\xi$ either increasing from $0$ or decreasing from 
	$\infty$), it provides the values for $V(\xi)$ and $C(\xi)$. 
	These in turn define the flow variables $u(t,x)$ and $c(t,x)$ along curves 
	$\xi=t^{-1/\lambda}x\equiv constant$, according to \eq{sim_var_u} and 
	\eq{sim_var_c}.
	The latter curves foliate the $(x,t)$-plane, with shapes depending on
	$\lambda$.
	\item The challenge is then to continue the trajectory $\Gamma_0$ so that 
	it connects $P_0$ to a critical point in the $(V,C)$-plane corresponding to 
	a vacuum state in physical space.
	\item There are two additional critical points $P_1=(-1,0)$ and 
	$P_2=(-\lambda,0)$ of the similarity ODEs \eq{V_ode}-\eq{C_ode} 
	located along $\{C=0\}$. According to \eq{sim_var_c} these points can 
	therefore correspond to a vacuum in the resulting Euler flow. 
	\item Our task is therefore to analyze how the trajectory $\Gamma_0$ can be
	connected to either $P_1$ or $P_2$. If $\Gamma_0$ itself joins 
	$P_0$ directly to $P_1$ or  $P_2$, then we obtain a globally defined and continuous
	Euler flow. If $\Gamma_0$ reaches neither $P_1$ nor 
	$P_2$, then we need to show that it is possible to jump, in an admissible manner, 
	from a point on $\Gamma_0$ to a point which is connected to either $P_1$ or 
	$P_2$ via a different trajectory $\Gamma_1$.
	The jump in the $(V,C)$-plane induces a corresponding discontinuity in the 
	physical flow variables, and ``admissible'' refers to the requirement that 
	the resulting discontinuity be an entropy admissible shock for the 
	Euler system \eq{mass}-\eq{mom}.
\end{itemize}
Several remarks are in order.
First, let us clarify which states in the $(V,C)$-plane can correspond to a 
vacuum in physical space. In addition to $P_1$ and $P_2$, the critical point 
$P_0$ at the origin is also located along $\{C=0\}$. However, we claim
that, for the solutions under consideration,  $P_0$ cannot 
serve as a vacuum state. To see this, recall that $P_0$ is
a star point for \eq{CV_ode} where every trajectory approaches tangent to
a straight line. The similarity ODE \eq{V_ode} gives that the leading
order behavior along a solution with $C\sim kV$ as $V\to0$ is
$V(\xi)\sim |\xi|^{-\lambda}$. In particular, approach to $P_0$ means
$|\xi|\to\infty$ if $\lambda>0$, or $|\xi|\to0$ if $\lambda<0$. Because of 
\eq{VC_at_t=0}, the only possibility for the trajectory 
$\Gamma_0$ (which emanates from $P_0$) to return to $P_0$ would 
be as $\xi\to-\infty$ if $\lambda>0$, or as $\xi\to0-$ if $\lambda<0$.
In particular, the resulting Euler flow would be defined at all points
$(t,x)$ with $t>0$.
However, in either case the sound speed $c(t,x)$, and hence the density 
$\rho(t,x)$ in the corresponding Euler flow (cf.\ \eq{sim_var_c}), would then
take strictly positive values at {\em all} points $(t,x)$ with $t>0$. 
With the initial vacuum region located along $\{x<0\}$, this behavior would imply 
infinite speed of propagation. This unphysical behavior shows that $P_0$ is not 
a candidate for a vacuum state. This leaves $P_1$ and $P_2$ as
the only potential vacuum states along the $V$-axis.

It turns out that there is one further possibility for reaching a vacuum state: 
In the particular cases $\text{Ma}=\pm\ell$, the trajectory $\Gamma_0$ 
lies along one of the straight lines $E_\pm:=\{C=\pm\ell^{-1}V\}$, which 
contain the critical points $P_3$ and $P_4$ (see Sections \ref{crit_pnt_lcns} 
and  \ref{P34}). While these points are 
located off the $V$-axis, the solutions along $E_\pm$ 
are such that $\xi C(\xi)\to0$ as $P_3$ and $P_4$ are approached. 
According to \eq{sim_var_c_alt}, this again yields 
approach to vacuum. 

With this we have identified the critical points $P_1$-$P_4$ that can serve as 
possible vacuum states. However, at this stage,
it is not obvious that it is possible to connect $P_0$ 
to one of them, continuously or not. Our analysis will show that this is indeed 
possible in case (I) (and this will require use of the last pair of
critical points $P_5$ and $P_6$ for \eq{CV_ode} as well). However, as noted in 
Section \ref{interpretation}, this is not always possible in case (II). In the latter
situation we obtain a ``solution'' of the original initial value problem which 
is only defined on a part of the $(x,t)$-plane.

Next, for case (I) we shall see that a subset of the relevant trajectories 
$\Gamma_0$ necessarily move off to infinity in the lower $(V,C)$ half-plane as 
$\xi\downarrow0$. An analysis of \eq{CV_ode} reveals that there is a one-to-one 
correspondence between such $\Gamma_0$-trajectories and asymptotic 
slopes $k\neq\pm1$ at infinity. To build complete Euler flows in these 
cases, $\Gamma_0$ must be continued into the upper half-plane 
coming in from infinity with the same asymptotic slope $k$ (cf.\ Section \ref{IAb}). 
This part of the analysis will be done in inverted 
coordinates $(V^{-1},C^{-1})$ and exploits a classic ODE 
result concerning approach toward non-simple equilibria.

We note that the outline above offers only two  scenarios: 
either the flow is continuous or it contains a single admissible shock. 
Furthermore, with the vacuum located to the left of the fluid, any 
admissible shock emanating from the initial interface must necessarily
be a 2-shock (i.e., with characteristics of the second family
running into it as time increases). This follows since no 1-characteristic emanates from
the vacuum interface. In contrast, 2-characteristics do emanate 
(tangentially) from the interface and can, under suitable conditions,
proceed to impinge on a propagating 2-shock. 

Thus, in addition to identifying and analyzing the critical points of
the similarity ODEs, we must also analyze the (self-similar) 
Rankine-Hugoniot relations. In Section \ref{ss_grps_rh_e_conds} we identify 
the possible locations in the $(V,C)$-plane of left and 
right states $P_\pm$ of admissible shocks. (For completeness we treat 
shocks of both families.) In order to argue for the appearance 
of admissible shocks we shall need to know how $P_+$, say, behaves 
as $P_-$ moves along certain trajectories of the similarity ODEs 
\eq{V_ode}-\eq{C_ode}; see Section \ref{hug_locus}.

Finally, having analyzed the critical points and jump relations for admissible 
shocks in self-similar flow, the resolution of the initial value problem 
\eq{u}, \eq{c}, \eq{vac_ss_data}
is reduced to the identification of suitable trajectories connecting 
the origin $P_0$ to one of the critical points $P_1$-$P_4$.
This is carried out in Sections \ref{I}-\ref{II} for 
cases (I) and (II). Since the phase portrait 
of \eq{CV_ode} changes at $\gamma=3$, we treat 
the cases $\gamma\gtrless 3$ separately. It turns out that case (I)
with $\gamma>3$ requires a further sub-division 
depending on the value of $\lambda$; see Section \ref{other_cases}.

\section{Critical points along the $V$-axis}\label{P0_P2}
The critical points are the points of intersection between the zero-levels 
\[\mathcal F:=\{(V,C)\,:\, F(V,C)=0\}\qquad\text{and}\qquad 
\mathcal G:=\{(V,C)\,:\, G(V,C)=0\}\]
of the functions defined in \eq{F} and \eq{G}, respectively. 
Note that $V=V_*$ (see \eq{V_*}) is a vertical asymptote for $\mathcal G$.
It turns out that there are up to seven points of intersection between 
$\mathcal F$ and $\mathcal G$, and we number these $P_i=(V_i,C_i)$, 
$i=0,\dots,6$. 

In this section we identify and analyze the critical points of \eq{CV_ode} 
along the $V$-axis. Since $F$ vanishes identically there, while 
$G(V,0)=0$ has the roots $V=0,-1,-\lambda$, there are always 
exactly three critical points located along the $V$-axis (recall the 
standing assumption that $\lambda\neq0,1$):
\[P_0:=(0,0) \qquad\qquad P_1:=(-1,0) \qquad\qquad P_2:=(-\lambda,0).\]  
As observed in Section \ref{outline}, all solutions of \eq{V_ode}-\eq{C_ode}
under consideration must tend to $P_0$ as $\xi\downarrow0$ 
or $\xi\uparrow\infty$, while $P_1$-$P_4$ 
provide possible end-states that describe approach to vacuum.

\subsection{Critical point $P_0=(0,0)$}\label{P0}
The linearization of \eq{CV_ode} at $P_0$ is $\frac{dC}{dV}=\frac C V$, 
so that $P_0$ is a star point (proper node) for all values of the parameters 
$\gamma$ and $\lambda$. According to Theorem 3.5 (iv) in \cite{ha},
the qualitative behavior of trajectories of the nonlinear ODE \eq{CV_ode}
near the origin agrees with that of its linearization there. Specifically, the 
slope at which a trajectory of \eq{CV_ode} approaches the origin uniquely 
determines the trajectory.


\subsection{Critical point $P_1=(-1,0)$}\label{P1}
Linearizing \eq{CV_ode} about $P_1$ yields 
\[\frac{dC}{dV}=-\ell^{-1}\frac{C}{1+V},\]
showing that $P_1$ is a saddle point. To obtain a more precise 
description we switch to coordinates 
\beq\label{P1_coords}
	(W,Z):=(1+V,C^2).
\eeq
A calculation shows that \eq{CV_ode} is transformed to
\beq\label{ZW_ode_1}
	\frac{dZ}{dW}=\frac{f_1(W,Z)}{g_1(W,Z)},
\eeq
where
\[f_1(W,Z)=2Z(Z-W^2+k_1W-k_0)\]
and
\[g_1(W,Z)=Z(W-W_*)-W(W-1)(W+\lambda-1),\]
with $W_*=1+V_*$ (cf.\ \eq{V_*} and \eq{ks}).
Linearizing \eq{ZW_ode_1} about $(W,Z)=(0,0)$ yields
\beq\label{P1_linzn}
	\frac{dZ}{dW}=\frac{a_1Z}{b_1W+c_1Z},
\eeq
where 
\[a_1=2\ell^{-1}(\lam-1)<0,\qquad
b_1=1-\lambda>0,\qquad
c_1=1+\ell(1-\lambda)>0.\]
The characteristic values of \eq{ZW_ode_1} are $a_1$ and $b_1$, 
with corresponding characteristic slopes
\[\sigma_1:=\frac{a_1-b_1}{c_1}
=-\frac{\gamma(\gamma-1)(1-\lambda)}{(\gamma-1)+2(1-\lambda)}
\qquad\text{and}\qquad \tau_1:=0,\]
respectively. Note that $\sigma_1<0$ (since $\lambda<1<\gamma$). Translating back to 
$(V,C)$-coordinates we therefore obtain the following:
$P_1=(-1,0)$ is a saddle point  for  \eq{CV_ode} at which two of the 
separatrices approach $P_1$ tangent to  
\beq\label{P1_1st_separatrix}
	C=\pm\sqrt{\sigma_1(1+V)}, \qquad\text{with $V<-1$,}
\eeq
and the other two separatrices lie along the $V$-axis to the left and right of 
$V=-1$. Focusing on the upper half-plane near $P_1$, let the separatrix there 
be denoted $\Sigma'$. We claim that $\Sigma'$ necessarily lies above
the level set $\{G=0\}$. Indeed, according to \eq{G} we have, to leading order, 
that $G$'s zero level through 
$P_1$ lies along $C=\sqrt{\sigma_{G}(1+V)}$, where $V<-1$ and 
$\sigma_{G}=(\lambda-1)/(1+\ell(1-\lambda))$.
A calculation shows that $|\sigma_1|>|\sigma_{G}|$, establishing the claim.
See Figure \ref{IA_fig_1g} for a schematic picture in Case (I) (i.e., $0<\lambda<1$
and $1<\gamma<3$). 

Next, consider a solution of the original similarity ODEs \eq{V_ode}-\eq{C_ode} 
moving along $\Sigma'$. Evaluating \eq{C_ode} with 
$C\sim\sqrt{\sigma_1(1+V)}$ we obtain that the leading order behavior 
along $\Sigma'$ near $P_1$ is given by
\beq\label{lob_P1}
	\frac{dC}{d\xi}\sim \frac{k_0}{\xi}\frac{1}{C},
\eeq
where $k_0$ is given in \eq{ks}.
It follows that the point $P_1$ is reached along $\Sigma'$ with a 
finite and non-zero $\xi$-value denoted $\xi_v$. Note that, 
as $\Sigma'$ is located in the upper half-plane, \eq{sign_cond}
implies that $\xi_v<0$.

The point $P_1$ will, for certain cases, provide the endpoint of 
the relevant solution to the similarity ODEs \eq{V_ode}-\eq{C_ode}.
In particular, since $C$ vanishes at $P_1$, \eq{sim_var_c} shows that
approaching $P_1$ with $\xi\to\xi_v$ corresponds to approaching a vacuum interface 
located along $x=\xi_v t^{1/\lambda}$ in the $(x,t)$-plane. 

Finally, consider the decay to vacuum in the corresponding Euler flow. 
Using \eq{sim_coord} and 
\eq{sim_var_c}, we calculate that for any fixed time $t>0$,
\beq\label{approach_to_vac_P1}
	(c^2(t,x))_x=2\lambda^{-2}t^{\frac{1}{\lambda}-2}
	\xi C(\xi)\left[C(\xi)+\xi C'(\xi)\right].
\eeq
Letting $x\downarrow \xi_v t^\frac{1}{\lambda}$, i.e., $\xi\downarrow\xi_v$, 
we get from \eq{lob_P1} and \eq{approach_to_vac_P1} that
\beq\label{approach_to_vac_P1_2}
	(c^2(t,x))_x\to 2\lambda^{-2}t^{\frac{1}{\lambda}-2}k_0\xi_v
	\qquad\text{as $x\downarrow \xi_v t^\frac{1}{\lambda}$.}
\eeq
Note that, as we restrict attention to $\lambda<1$
(so that $k_0<0$) and since $\xi_v<0$, \eq{approach_to_vac_P1_2} 
gives a positive value for $(c^2(t,x))_x$ along the interface,
as must be the case with the vacuum located to the left.

Thus, whenever the approach to vacuum in self-similar Euler flow
corresponds to approaching $P_1$ along the separatrix $\Sigma'$, 
then the decay to vacuum from within the fluid is given by a 
physical singularity: $c^2$ is Lipschitz continuous with respect to 
$x$ at each fixed time $t>0$.

\subsection{Critical point $P_2=(-\lambda,0)$}\label{P2}
Linearizing \eq{CV_ode} about $P_2$ yields 
\beq\label{P2_linzn}
	\frac{dC}{dV}=\ell^{-1}\frac{C}{\lambda+V},
\eeq
showing that $P_2$ is a nodal point. For $\gamma>3$ (i.e., $\ell^{-1}<1$) 
two trajectories approach $P_2$ along $\{C=0\}$ and all other trajectories 
approach $P_2$ vertically. For $1<\gamma<3$ the opposite holds, while 
$P_2$ is a star point (proper node) when $\gamma=3$.

\section{Critical points off the $V$-axis}\label{P3_P6}
The remaining critical points $P_3$-$P_6$ may be obtained 
by solving $G(V,C)=0$ for $C^2$ in terms of $V$, and substituting 
the result into the equation $F(V,C)=0$; this yields a quadratic
polynomial in $V$ (see below). According to \eq{symms} the critical 
points $P_3$-$P_6$ come in pairs located symmetrically about the 
$V$-axis. With $P_3$ and $P_5$ denoting the ones in the upper 
half-plane, we have
\[P_3=(V_3,C_3),\qquad P_4=(V_3,-C_3),\qquad P_5=(V_5,C_5),
\qquad P_6=(V_5,-C_5),\]
where $C_3,C_5>0$. 
We verify below that $P_3$ and $P_4$ are present 
for all values of $\lambda$ and $\gamma$, while $P_5$ and $P_6$ 
are present whenever $\gamma\neq 3$.

Restricting attention to $P_3$ and $P_5$, we proceed to determine 
these. From $G(V,C)=0$ we have 
\beq\label{C^2_G}
	C^2=\textstyle\frac{V(1+V)(\lambda+V)}{V-V_*}=:g(V).
\eeq
Substituting \eq{C^2_G} into $F(V,C)=0$, and recalling that we now seek 
critical points off the $V$-axis, we obtain a quadratic equation for $V$.
For $\gamma\neq 3$ its roots are given by
\beq\label{V_3}
	V_3=-\textstyle\frac{2\lambda}{\gamma+1}=-\frac{\lambda\ell}{\ell+1},
\eeq
and
\beq\label{V_5}
	V_5=\textstyle\frac{2}{\gamma-3}=\frac{\ell}{1-\ell}.
\eeq
For $\gamma=3$ the quadratic equation degenerates to a linear equation with the 
single root $V_3|_{\gamma=3}=-\frac{\lambda}{2}.$

It remains to verify that these roots satisfy $g(V_3),g(V_5)\geq 0$ (cf.\ \eq{C^2_G}),
so that $P_3$ and $P_5$ are present. A calculation shows that 
\[g(V_3)=\textstyle\frac{\lambda^2}{(\ell+1)^2}\geq 0\qquad\text{and}\qquad
g(V_5)=\textstyle\frac{1}{(\ell-1)^2}>0.\]
It follows that $P_3$ and $P_4$ are always present with
\beq\label{P_3_4}
	P_3=(-\textstyle\frac{\lambda\ell}{\ell+1},\textstyle\frac{|\lambda|}{\ell+1}),\qquad 
	P_4=(-\textstyle\frac{\lambda\ell}{\ell+1},-\textstyle\frac{|\lambda|}{\ell+1}),
\eeq
while $P_5$ and $P_6$ are present if and only if $\gamma\neq 3$ ($\ell\neq 1$), and 
\beq\label{P_5_6}
	P_5=\textstyle(\frac{\ell}{1-\ell},\textstyle\frac{1}{|\ell-1|}),\qquad
	P_6=\textstyle(\frac{\ell}{1-\ell},-\textstyle\frac{1}{|\ell-1|})\qquad (\gamma\neq 3).
\eeq

%
%
%
%
%

\subsection{Explicit trajectories and locations of $P_3$-$P_6$}\label{crit_pnt_lcns}
A direct calculation shows that the reduced similarity ODE \eq{CV_ode}
always admits exactly two straight-line trajectories
\[E_\pm:=\{(V,C)\,|\, C=\pm\ell^{-1}V\}.\]
These will be useful in delimiting solution behaviors. 
It is immediate to verify that the critical points $P_3$-$P_6$ off the 
$V$-axis are always located on $E_\pm$, and that
\begin{itemize}
	\item $P_3\in E_-\Leftrightarrow\lambda>0$;
	\item $P_3\in E_+\Leftrightarrow\lambda<0$;
	\item $P_5\in E_-\Leftrightarrow\gamma<3$;
	\item $P_5\in E_+\Leftrightarrow\gamma>3$.
\end{itemize}
Also, by symmetry we have 
\[P_4\in E_\pm\Leftrightarrow P_3\in E_\mp ,\qquad\text{and}\qquad 
P_6\in E_\pm\Leftrightarrow P_5\in E_\mp.\]
For later reference we record the ODE for $V(\xi)$ 
when evaluated along $E_\pm$: with $C=\pm\ell^{-1}V$ equation
\eq{V_ode} becomes
\beq\label{V_ode_along_E_+}
	\frac{dV}{d\xi}=-\frac{1}{\xi}\frac{V(V-V_3)}{(V+\frac{\ell}{1+\ell})}.
\eeq
 
We also record the locations of $P_3$-$P_6$ 
relative to the critical lines $L_\pm=\{C=\pm(1+V)\}$. 
First, \eq{P_5_6} shows that whenever 
$P_5$ and $P_6$ are present (i.e., when $\gamma\neq 3$), 
these points are necessarily located on $L_-\cup L_+$.
To describe the locations of $P_3$ and $P_4$ relative to $L_\pm$ 
we  introduce the open cone
\[\mathcal K:=\{(V,C)\,|\, C^2<(1+V)^2\},\]
whose boundary is $L_-\cup L_+$. Using \eq{P_3_4} we calculate 
that 
\beq\label{P_3_in_K}
	P_3,P_4\in \mathcal K\qquad\Longleftrightarrow\qquad
	(\lambda-1)[(\gamma-3)\lambda+(\gamma+1)]<0.
\eeq

\subsection{Critical points $P_3$ and $P_4$}\label{P34}
Due to the symmetries \eq{symms} it suffices to analyze $P_3$.
In this subsection, unless indicated differently, all quantities are 
evaluated at $P_3$, and the subscript `$3$' is mostly suppressed. 
Linearizing \eq{CV_ode} about $P_3$ yields 
\beq\label{P3_linz'n}
	\frac{dc}{dv}=\frac{F_Vv+F_Cc}{G_Vv+G_Cc},
\eeq
where $v:=V-V_3$, $c:=C-C_3$, and the partial derivatives 
(evaluated at $P_3$) are given by 
\begin{align}
	F_C&=2C^2\label{F_C_3}\\
	F_V&=C(k_1-2(1+V))\label{F_V_3}\\
	G_C&=2C(V-V_*)\label{G_C_3}\\
	G_V&=C^2-(3V^2+2(1+\lambda)V+\lambda).\label{G_V_3}
\end{align}
In general, the Wronskian $W$ and discriminant $R^2$ 
of \eq{P3_linz'n} are defined by
\beq\label{W_R}
	W:=F_CG_V-F_VG_C \qquad\text{and}\qquad 
	R^2:= (F_C+G_V)^2-4W.
\eeq
A direct calculation shows 
that $W$ at $P_3$ takes the value
\[W=\ell C^2(\lambda-1)
[(\gamma-3)\lambda+(\gamma+1)].\]
In particular, we get that the 
Wronskian is strictly negative whenever $\lambda<1$ and $1<\gamma<3$. 
This implies that $P_3$ is necessarily a saddle point in these cases. 

%

\subsection{Critical points $P_5$ and $P_6$}\label{P56}
It turns out that, for our needs in resolving self-similar Euler flows,
the point $P_5=(V_5,C_5)$ is relevant only when
$0<\lambda<1$ and $1<\gamma<3$, while $P_6=(V_5,-C_5)$
is relevant only when $0<\lambda<1$, $\gamma>3$. 
While most of the calculations for the former case carry over to the latter,
there are also differences due to the different ranges for $\gamma$.
In particular, the phase portrait of \eq{CV_ode} changes
at $\gamma=3$, and for the second case we shall need to consider two 
further sub-cases.

\subsubsection{$P_5$ with $0<\lambda<1$, $1<\gamma<3$}\label{P5}
In this case $P_5=(\frac{2}{\gamma-3},\frac{\gamma-1}{3-\gamma})$ 
is located in the second quadrant and to the left of $V=-1$.
Also, $k_1>0$, $k_0<0$, and $V_*>0$, cf.\ \eq{V_*}-\eq{ks}. Below,
unless indicated differently, all quantities are evaluated at $P_5$, and
the subscript `$5$' is suppressed in most of the expressions. 

To determine the behavior of \eq{CV_ode} near $P_5$ we need the signs of various
quantities given in terms of the partial derivatives of $F$ and $G$.
Recall that $P_5$ is a ``triple point'' which belongs to $L_-$, $\{F=0\}$, and $\{G=0\}$; 
in particular, at $P_5$ we have
\begin{align}
	C&=-(1+V)\label{C_V_5_1}\\
	C^2&=(1+V)^2-k_1(1+V)+k_0\label{C_V_5_2}\\
	C^2&=\textstyle\frac{V(1+V)(\lambda+V)}{V-V_*}.\label{C_V_5_3}
\end{align}
At $P_5$ we then have (using \eq{C_V_5_1})
\begin{align}
	F_C&=2C^2>0\label{F_C_5}\\
	F_V&=C(k_1-2(1+V))\equiv C(k_1+2C)>0\label{F_V_5}\\
	G_C&=2C(V-V_*)\equiv -2V(\lambda+V)<0\label{G_C_5}\\
	G_V&=C(\lambda-1+2V-V_*)<0.\label{G_V_5}
\end{align}
Here, $F_C$ and $F_V$ are calculated from \eq{F}, 
while $G_C$ and $G_V$ are calculated from \eq{G} using 
\[C(V_*-V)=V(\lambda+V),\]
which is a consequence of \eq{C_V_5_1} and \eq{C_V_5_3}. 
We next determine the relative positions of the curves $L_-$, $\{F=0\}$, $\{G=0\}$,
as well as the straight-line trajectory $E_-=\{C=-\ell^{-1}V\}$, near $P_5$.
(See Figure \ref{IA_fig_1g}.) 
\begin{lemma}\label{rel_posns_5}
	When $0<\lambda<1$, $1<\gamma<3$ we have the following relations at $P_5$:
	\beq\label{slopes_5}
		\textstyle-\frac{F_V}{F_C}<-1<-\frac{G_V}{G_C}<-\ell^{-1}.
	\eeq
\end{lemma}
\begin{proof}
	Differentiating the relation $\ell^{-1}G(V,-(1+V))=F(V,-(1+V))$
	(cf.\ \eq{non_obvious_reln}) with respect to $V$, 
	and applying \eq{F_C_5}-\eq{F_V_5}, give
	$\ell^{-1}(G_V-G_C)=F_V-F_C=k_1C>0$. Using the signs in 
	\eq{F_C_5}-\eq{G_V_5} then gives the two leftmost inequalities in 
	\eq{slopes_5}. For the rightmost inequality, a direct calculation 
	using \eq{G_C_5}-\eq{G_V_5} shows that 
	it reduces to the inequality $0<2+(\ell-1)(1-\lambda)$, which is
	satisfied in the case under consideration.
\end{proof}
Next, consider the linearization of \eq{CV_ode} at $P_5$,
viz.\ \eq{P3_linz'n} with $v:=V-V_5$, $c:=C-C_5$, and partial derivatives
given by \eq{F_C_5}-\eq{G_V_5}. The Wronskian $W$ and discriminant $R^2$ 
are given in \eq{W_R}.
Note that $W>0$ in the present case according to \eq{F_C_5}-\eq{G_V_5} 
and Lemma \ref{rel_posns_5}. Also, a direct calculation shows that the 
discriminant is given by
\beq\label{discrmn}
	R^2=C^2(2+(\ell-3)(1-\lambda))^2.
\eeq
With $0<\lambda<1$ and $1<\gamma<3$ we have $2+(\ell-3)(1-\lambda)>0$, so that $R^2>0$. 
With $R:=+\sqrt{R^2}>0$ we set
\beq\label{L_12}
	L_{1,2}=\textstyle\frac{1}{2G_C}(F_C-G_V\pm R)
\eeq
and
\beq\label{E_12}
	E_{1,2}=\textstyle\frac{1}{2G_C}(F_C+G_V\pm R),
\eeq
and chose signs so that 
\beq\label{Es}
	|E_1|<|E_2|. 
\eeq
Note that the signs $\pm$ in \eq{L_12} and
in \eq{E_12} agree; $L_1$ and $L_2$ are referred to as
the {\em primary} and {\em secondary} slopes (or directions), respectively. 
In terms of these we have that trajectories of \eq{CV_ode} near $P_5$ 
approach one of the curves
\[(c-L_1v)^{E_1}=\text{constant}\times (c-L_2v)^{E_2}.\]

Since $W,R>0$, $P_5$ is a node: all trajectories of \eq{CV_ode} approaching 
$P_5$ do so with slope equal to the primary slope $L_1$, except two 
which approach $P_5$ with slope $L_2$. We proceed to calculate 
the primary and secondary slopes. (Since the special straight-line 
trajectory $E_-$ passes through $P_5$, we already know that one of these 
slopes must be given by $-\ell^{-1}$; we verify below that this is the secondary 
slope at $P_5$.) First, from \eq{F_C_5} and \eq{G_V_5} we have
\[F_C+G_V=-C(1-\lambda+V_*+2)<0,\]
so that
\[E_{1,2}=\textstyle\frac{1}{2G_C}(-|F_C+G_V|\pm \sqrt{|F_C+G_V|^2-4W}).\]
Since $W>0$, the choice in \eq{Es} requires that subscript `1' corresponds to 
the `$+$' sign. Thus,
\beq\label{prmr_scdr_slopes_5_1}
	L_1=\textstyle\frac{1}{2G_C}(F_C-G_V+R)
	\qquad\text{and}\qquad
	L_2=\textstyle\frac{1}{2G_C}(F_C-G_V-R).
\eeq
By evaluating $R$ and the various partial derivatives at $P_5$ we obtain
\beq\label{prmr_scdr_slopes_5_2}
	L_1=-\textstyle\frac{2\ell+(\ell-1)^2(1-\lambda)}{2\ell[1+(\ell-1)(1-\lambda)]}
	\qquad\text{and}\qquad
	L_2=-\ell^{-1}.
\eeq
Also, direct calculations verify the following inequalities (which include those in 
Lemma \ref{rel_posns_5}):
 \beq\label{prmr_scdr_slopes_5_3}
	\textstyle-\frac{F_V}{F_C}<-1<-\textstyle\frac{G_V}{G_C}<L_1<-\ell^{-1}\equiv L_2.
\eeq
(See Figure \ref{IA_fig_1g} where the trajectories $I\!I\!I$-$I\!V$ and $\Sigma'$ approach
$P_5$ with slope $L_1$).

Finally, we note that whenever $P_5$ is approached by a solution of the 
similarity ODEs \eq{V_ode}-\eq{C_ode}, the independent variable $\xi$ tends 
to a finite value. This follows since \eq{V_ode}, say, evaluated along 
$C-C_5=L(V-V_5)$, with $L=L_1$ or $L=L_2$, gives
\[\textstyle\frac{dV}{d\xi}\approx\frac{A}{\xi}\qquad\text{for $V\approx V_5$,}\] 
where $A=(G_V+LG_C)/2(1+L)(1+V_5)$ is a finite constant.

\subsubsection{$P_6$ with $0<\lambda<1$, $\gamma>3$}\label{P6}
In this case $P_6$ (located in the lower half-plane) is given by the same expression 
as $P_5$ above, i.e., $P_6=(\textstyle\frac{2}{\gamma-3},\frac{\gamma-1}{3-\gamma}).$
$P_6$ is thus located in the fourth quadrant and to the right of $V=V_*$.
We now have $k_1<0$, $k_0<0$, and $V_*>0$, cf.\ \eq{V_*}-\eq{ks}.
Unless indicated differently, all quantities in this subsection are evaluated 
at $P_6$, and the subscript `$6$' is mostly suppressed.

Since the expression for $P_6$ coincides with that of $P_5$  in Section \ref{P5}, 
the calculations there apply verbatim. In particular, $P_6$ is 
a triple point belonging to $L_-\cap\{F=0\}\cap\{G=0\}$, and the partials of 
$F$ and $G$ are again given by the expressions in \eq{F_C_5}-\eq{G_V_5}. 
Also, the signs of $F_C$, $F_V$, and $G_C$ agree with those displayed 
in \eq{F_C_5}-\eq{G_C_5}: $F_C>0$, $F_V>0$, and $G_C<0$. 
However, $G_V=C(\lambda-1+2V-V_*)$ may now be of either sign depending 
on $\lambda$ and $\gamma$. 

We next record the relative positions of the curves $L_-$, $\{F=0\}$, $\{G=0\}$,
as well as the straight-line trajectory $E_-=\{C=-\ell^{-1}V\}$, near $P_6$.
(The proof is similar to that of Lemma \ref{rel_posns_5}.)
\begin{lemma}\label{rel_posns_6}
	When $0<\lambda<1$, $\gamma>3$ we have the following relations at $P_6$:
	\beq\label{slopes_6}
		\textstyle-\ell^{-1}<-\frac{F_V}{F_C}<-1<-\frac{G_V}{G_C}.
	\eeq
\end{lemma}
Next, consider the linearization of \eq{CV_ode} at $P_6$, 
viz.\ \eq{P3_linz'n} with $v:=V-V_6$, $c:=C-C_6$, and partial derivatives
given by \eq{F_C_5}-\eq{G_V_5}. The Wronskian $W$ and discriminant $R^2$ 
are defined as in \eq{W_R}.
Again, from Lemma \ref{rel_posns_6} and the signs
$F_C>0$ and $G_C<0$, we obtain the Wronskian
\[W=F_CG_V-F_VG_C>0,\]
while the discriminant is given by $R^2=C^2(2+(\ell-3)(1-\lambda))^2$. 
Thus, $R^2\geq0$, so that
\beq\label{R_IB}
	R=-C|2+(\ell-3)(1-\lambda)|,
\eeq
where $C=C_6<0$.
Disregarding the particular case $R=0$ for now, we consider two sub-cases 
depending on the sign of $2+(\ell-3)(1-\lambda)$. For $\gamma>3$ fixed we set
\beq\label{lambda_hat}
	\hat\lambda:=\frac{\gamma-3}{3\gamma-5},
\eeq
and we have
\begin{enumerate}
	\item[(i)] If $\hat\lambda<\lambda<1$, then 
	$R=-C(2+(\ell-3)(1-\lambda))>0$.
	\item[(ii)] If $0<\lambda<\hat\lambda$, then 
	$R=C(2+(\ell-3)(1-\lambda))>0$.
\end{enumerate}
(Note that $\gamma>3$ implies $0<\hat\lambda<1$, so that 
both cases (i) and (ii)  are possible.)
Arguing as for $P_5$ above we get that $P_6$ is a node and that the 
following holds.
\begin{itemize}
	\item In case (i), the primary and secondary slopes at $P_6$  are given by
	\beq\label{prmr_scdr_slopes_6_1}
		L_1=-\textstyle\frac{2\ell+(\ell-1)^2(1-\lambda)}{2\ell[1+(\ell-1)(1-\lambda)]}
		\qquad\text{and}\qquad L_2=-\ell^{-1},\quad\text{respectively.}
	\eeq
	\item In case (ii), the primary and secondary slopes at $P_6$  are given by
	\beq\label{prmr_scdr_slopes_6_2}
		L_1=-\ell^{-1}
		\qquad\text{and}\qquad L_2=-
		\textstyle\frac{2\ell+(\ell-1)^2(1-\lambda)}{2\ell[1+(\ell-1)(1-\lambda)]},
		\quad\text{respectively.}
	\eeq
\end{itemize}
In the special case that $R=0$ we have that $P_6$ is a degenerate node for
which the primary and secondary slopes take the same value $-\ell^{-1}$.
Next, direct calculations verify the following inequalities (which include those in 
Lemma \ref{rel_posns_6}).
\begin{itemize}
	\item For case (i):
	\beq\label{prmr_scdr_slopes_6_3}
		-\ell^{-1}<L_1<-\textstyle\frac{F_V}{F_C}<-1<-\textstyle\frac{G_V}{G_C}
	\eeq
	\item For case (ii):
	\beq\label{prmr_scdr_slopes_6_4}
		L_2<-\ell^{-1}<-\textstyle\frac{F_V}{F_C}<-1<-\textstyle\frac{G_V}{G_C}.
	\eeq
\end{itemize}
Finally, arguing as for $P_5$ above, one may verify the following: 
Whenever $P_6$ is approached by a solution of the 
similarity ODEs \eq{V_ode}-\eq{C_ode}, the independent variable $\xi$ tends 
to a finite value.

\section{Jump and entropy conditions}\label{ss_grps_rh_e_conds}
This section considers the jump relations and entropy conditions for shocks in
self-similar solutions to \eq{mass}-\eq{mom}. We first characterize the regions 
in the $(V,C)$-plane that can be connected by admissible 1-shocks and 2-shocks
We then analyze the behavior of  the ``Hugoniot locus'' of a  
trajectory of \eq{CV_ode}, cf.\ Definition \ref{hug_loc_def}.

\subsection{Jump relations in self-similar variables}\label{rh_e_conds}
First let $(\rho,u)$ be a general solution of \eq{mass}-\eq{mom} in which a 
discontinuity propagates along $x=\mathcal X(t)$. 
The Rankine-Hugoniot jump relations are 
\beq\label{rh}	
	\dot{\mathcal X}\lj\rho\rj = \lj \rho u\rj \qquad\text{and}\qquad
	\dot{\mathcal X}\lj\rho u\rj = \lj \rho u^2+a^2\rho^\gamma\rj, 
\eeq
where we use the convention that for any quantity $q=q(t,x)$,
\[\lj q\rj:=q_+-q_-\equiv q(t,\mathcal X(t)+)-q(t,\mathcal X(t)-).\]
In what follows, a subscript `$-$' (`$+$') always refers to evaluation 
at the left (right) of a discontinuity in physical space.
The entropy conditions then require that 
\beq\label{e}
	u_--c_-> \dot{\mathcal X}> u_+-c_+\quad\text{(1-shock)}\qquad\text{ and}\qquad
	u_-+c_-> \dot{\mathcal X}> u_++c_+\quad\text{(2-shock).}
\eeq
Now assume that the solution under consideration is a self-similar solution of the 
form \eq{sim_var_u}-\eq{sim_var_c}, and that the 
shock propagates along the path $\xi\equiv \bar \xi$. Then
$\mathcal X(t)=\bar\xi t^\frac{1}{\lambda}$ and 
$\dot{\mathcal X}=\textstyle\frac{\mathcal X}{\lambda t}$.
Expressing the density $\rho$ in terms of the sound speed $c$ yields
$\rho=(a\sqrt{\gamma})^{-\ell} c^\ell$ ,
and using these relations in \eq{rh} we obtain the Rankine-Hugoniot 
relations \eq{rh} in $(V,C)$-variables: 
\beq\label{rh_VC}	
	\lj|C|^\ell(1+V)\rj = 0 \qquad\text{and}\qquad
	\lj|C|^\ell((1+V)^2+\textstyle\frac{1}{\gamma}C^2)\rj = 0,
\eeq
where $[\![\cdot]\!]$ now denotes jump across $\xi=\bar \xi$. Setting
\beq\label{RWM}
	R:=|C|^\ell,\qquad W:=1+V,\qquad M:=RW,
\eeq
these jumps relations are equivalent to
\beq\label{rh_RWM}	
	\lj M\rj = 0 \qquad\text{and}\qquad
	\lj\textstyle\frac{M^2}{R}+\textstyle\frac{1}{\gamma}R^\gamma\rj = 0.
\eeq
We record the following:
\begin{lemma}\label{f_m}
	For a given constant $m$, define the function $f_m:(0,\infty)\to(0,\infty)$ by
	\[f_m(R):=\textstyle\frac{m^2}{R}+\textstyle\frac{1}{\gamma}R^\gamma,\]
	and set
	\[R^*:=|m|^\frac{2}{\gamma+1}.\]
	Then $f_m$ is decreasing on $(0,R^*)$, increasing on $(R^*,\infty)$, and
	\[\lim_{R\to0+}f_m(R)=\lim_{R\to\infty}f_m(R)=+\infty.\]
\end{lemma}
As we consider flows for $t>0$, the entropy conditions \eq{e} for shocks 
propagating in a self-similar solution along $\xi\equiv\bar\xi$, take the form
\beq\label{1_shock}
	\textstyle\frac{\bar \xi}{\lambda}(C_--V_-)>\frac{\bar \xi}{\lambda}
	>\textstyle\frac{\bar \xi}{\lambda}(C_+-V_+)\qquad\text{for a 1-shock, and}
\eeq
\beq\label{2_shock}
	-\textstyle\frac{\bar \xi}{\lambda}(C_-+V_-)>\frac{\bar \xi}{\lambda}
	>-\textstyle\frac{\bar \xi}{\lambda}(C_++V_+)\qquad\text{for a 2-shock.}
\eeq
\begin{definition}
	Let $\lambda$ and $\bar\xi$ be fixed. We say that the pair of points 
	$P_-=(V_-,C_-)$ and $P_+=(V_+,C_+)$ defines an {\em admissible 
	self-similar 1-shock} with similarity parameter $\lambda$, left state $P_-$, 
	right state $P_+$, and propagating along $x=\bar\xi t^\frac{1}{\lambda}$, 
	provided \eq{rh_VC} and \eq{1_shock} are met. Admissible self-similar 2-shocks
	are defined similarly by requiring that \eq{rh_VC} and \eq{2_shock} are met.
\end{definition}
Next, we consider the locations of possible right and left states for admissible 
self-similar 1- and 2-shocks. We break down this issue into four cases depending
on $\frac{\bar \xi}{\lambda}\gtrless0$ and the type of the shock.

\subsubsection{Admissible self-similar shocks with $\frac{\bar \xi}{\lambda}>0$}\label{xi/lambda>0}
\begin{figure}[ht]
	\centering
	\includegraphics[width=8cm,height=8cm]{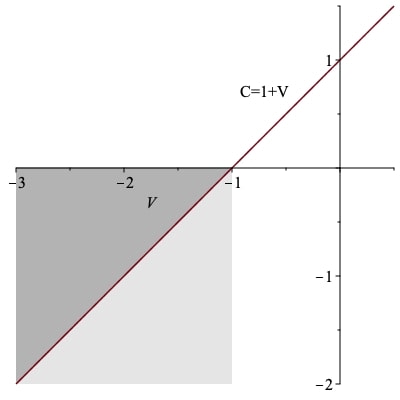}
	\caption{The regions $S_-^1=\{(V,C): 1+V<C<0\}$ (dark grey) and
	$S_+^1=\{(V,C): C<1+V<0\}$ (light grey) of 
	possible left and right states, respectively, of an admissible 
	1-shock when $\frac{\bar \xi}{\lambda}>0$.}\label{S1_fig}
\end{figure}

According to \eq{sign_cond} we are now only considering points in the lower half-plane
$\{C\leq0\}$. 
Let $\lambda$ and $\bar\xi$ be fixed with $\frac{\bar \xi}{\lambda}>0$, and 
define the regions
\[S^1_-:=\{(V,C): 1+V<C<0\},\qquad
S^1_+:=\{(V,C): C<1+V<0\};\]
see Figure \ref{S1_fig}. 
Note that, as $\frac{\bar \xi}{\lambda}>0$,
the points in $S^1_-$ ($S^1_+$, respectively) satisfy the leftmost (rightmost, respectively)
inequality in the entropy condition \eq{1_shock} for a 1-shock. The following proposition shows that 
these regions provides the possible locations for left and right states for
admissible self-similar 1-shocks in this case.
\begin{proposition}[Admissible self-similar 1-shocks with $\frac{\bar \xi}{\lambda}>0$]
\label{SS_1_shocks_with_xi/lambda>0}
	Assume $\frac{\bar \xi}{\lambda}>0$; then:
	\begin{enumerate}
	\item Whenever $P_-\in S_-^1$ there is a unique $P_+\in S^1_+$ 
		such that the pair $(P_-,P_+)$ satisfies the Rankine-Hugoniot relations \eq{rh_VC}.
	\item Conversely, whenever $P_+\in S^1_+$  there is a unique $P_-\in S^1_-$ 
		such that the pair $(P_-,P_+)$ satisfies the Rankine-Hugoniot relations \eq{rh_VC}.
	\end{enumerate}
	In either case, $(P_-,P_+)$  is an admissible 
	self-similar 1-shock with similarity parameter $\lambda$, left state $P_-$, 
	right state $P_+$, and propagating along $x=\bar\xi t^\frac{1}{\lambda}$.
\end{proposition}
\begin{proof}
	Let $P_-\in S_-^1$, i.e., $W_-<C_-<0$. It is convenient to employ the variables 
	$R$, $W$, $M$ defined in \eq{RWM}. Note that the admissibility conditions 
	\eq{1_shock} for 1-shocks take the form
	\beq\label{xi/lambda>0_1_shock}
		C_->W_-\qquad\text{and}\qquad C_+<W_+\qquad\qquad\text{(1-shock).}
	\eeq
	Also, for part (1) we have by assumption that
	\beq\label{assump1}
		M_-=R_-W_-<0,
	\eeq
	and that the Rankine-Hugoniot conditions \eq{rh_RWM} amount to the identities
	\beq\label{rh_1^-}
		f_-(R_+)=f_-(R_-)\qquad\text{and}\qquad R_+W_+=M_-,
	\eeq
	where, in the notation of Lemma \ref{f_m}, $f_-=f_{M_-}$.
	To establish part (1) we shall to argue that there is a unique
	unique value $R_+$ satisfying \eq{rh_1^-}${}_1$, and that, upon setting
	\beq\label{1^-_C_W}
		C_+:=-R_+^\frac{1}{\ell}\qquad\text{and}\qquad W_+:=\textstyle\frac{M_-}{R_+},
	\eeq
	(in particular, so that \eq{rh_1^-}${}_2$ is met) we have $P_+=(V_+,C_+)\in S^1_+$. 
	To show this, we first observe that 
	\[f_-'(R_-)=|C_-|^2-|W_-|^2<0,\]
	since $W_-<C_-<0$. It follows from Lemma \ref{f_m} that 
	\beq\label{R_-<R^*_}
		R_-<R^*_-:=|M_-|^\frac{2}{\gamma+1},
	\eeq
	and that there is a unique $R_+>R^*_-$ satisfying $f_-(R_+)=f_-(R_-)$. Defining 
	$C_+$ and $W_+$ according to \eq{1^-_C_W}, it remains to argue that
	$W_+<0$ and that $W_+>C_+$. The former inequality is immediate from 
	\eq{1^-_C_W}${}_2$ since $M_-<0$ and $R_+>0$; for the latter we have, according to \eq{1^-_C_W},
	that
	\[W_+>C_+\quad\Leftrightarrow\quad M_->-R_+^\frac{\gamma+1}{2}
	\quad\Leftrightarrow\quad M_-^2<R_+^{\gamma+1}\quad\Leftrightarrow\quad R^*_-<R_+,\]
	which was established just above. Thus, \eq{xi/lambda>0_1_shock} holds by 
	construction, finishing the proof of part (1).
	
	The proof of part (2) is similar: defining $f_+:=f_{M_+}$, we have $f_+'(R_+)>0$, so that 
	Lemma \ref{f_m} implies the existence of a unique value $R_-$ satisfying 
	$R_-<R^*_+:=|M_+|^\frac{2}{\gamma+1}<R_+$ so that \eq{rh_1^-}${}_1$ is met.
	Defining $C_-:=-R_-^\frac{1}{\ell}$ and $W_-:=\frac{M_+}{R_-}$, and arguing as above, establishes 
	part (2).
\end{proof}

\begin{figure}[ht]
	\centering
	\includegraphics[width=8cm,height=8cm]{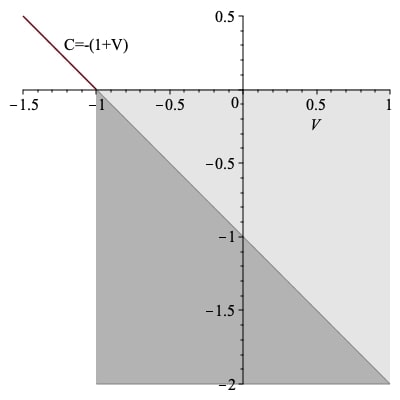}
	\caption{The regions $S^2_-:=\{(V,C): C<-(1+V)<0\}$ (dark grey) 
	and $S^2_+:=\{(V,C): -(1+V)<C<0\}$ (light grey) of 
	possible left and right states, respectively, of an admissible 2-shock when $\frac{\bar \xi}{\lambda}>0$.}\label{S2_fig}
\end{figure}

For the corresponding statement for 2-shocks we define the regions
\[S^2_-:=\{(V,C): C<-(1+V)<0\},\qquad
S^2_+:=\{(V,C): -(1+V)<C<0\};\]
see Figure \ref{S2_fig}.

As above, since $\frac{\bar \xi}{\lambda}>0$,
the points in $S^2_-$ ($S^2_+$, respectively) satisfy the leftmost (rightmost, respectively)
inequality in the entropy condition \eq{2_shock} for a 2-shock. Arguing as in the proof of Proposition
\ref{SS_1_shocks_with_xi/lambda>0}, we have the following result.
\begin{proposition}[Admissible self-similar 2-shocks with $\frac{\bar \xi}{\lambda}>0$]
\label{SS_2_shocks_with_xi/lambda>0}
	Assume $\frac{\bar \xi}{\lambda}>0$; then:
	\begin{enumerate}
	\item Whenever $P_-\in S_-^2$ there is a unique $P_+\in S^2_+$ 
		such that the pair $(P_-,P_+)$ satisfies the Rankine-Hugoniot relations \eq{rh_VC}.
	\item Conversely, whenever $P_+\in S^2_+$  there is a unique $P_-\in S^2_-$ 
		such that the pair $(P_-,P_+)$ satisfies the Rankine-Hugoniot relations \eq{rh_VC}.
	\end{enumerate}
	In either case, $(P_-,P_+)$  is an admissible 
	self-similar 2-shock with similarity parameter $\lambda$, left state $P_-$, 
	right state $P_+$, and propagating along $x=\bar\xi t^\frac{1}{\lambda}$.
\end{proposition}

\subsubsection{Admissible self-similar shocks with $\frac{\bar \xi}{\lambda}<0$}\label{xi/lambda<0}
\begin{figure}[ht]
	\centering
	\includegraphics[width=8cm,height=8cm]{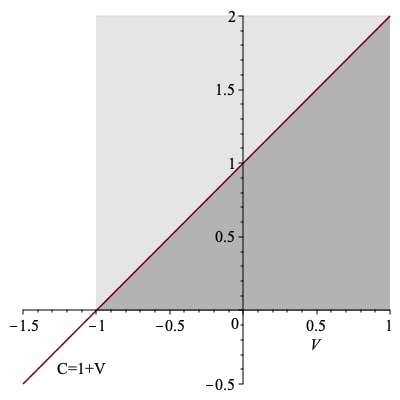}
	\caption{The regions $T^1_-=\{(V,C): 0<C<1+V\}$ (dark grey) 
	and $T^1_+:=\{(V,C): 0<1+V<C\}$ (light grey) of 
	possible left and right states, respectively, of an admissible 
	1-shock when $\frac{\bar \xi}{\lambda}<0$.}\label{T1_fig}
\end{figure}

According to \eq{sign_cond} we are now only considering points in the upper half-plane
$\{C\geq0\}$. Let $\lambda$ and $\bar\xi$ be fixed and such that $\frac{\bar \xi}{\lambda}<0$, and 
define the regions
\[T^1_-:=\{(V,C): 0<C<1+V\},\qquad
T^1_+:=\{(V,C): 0<1+V<C\};\]
see Figure \ref{T1_fig}. 
Since $\frac{\bar \xi}{\lambda}<0$,
the points in $T^1_-$ ($T^1_+$, respectively) satisfy the leftmost (rightmost, respectively)
inequality in the entropy condition \eq{1_shock} for a 1-shock. The proof of the following proposition
is similar to that of Proposition \ref{SS_1_shocks_with_xi/lambda>0}.
\begin{proposition}[Admissible self-similar 1-shocks with $\frac{\bar \xi}{\lambda}<0$]
\label{SS_1_shocks_with_xi/lambda<0}
	Assume $\frac{\bar \xi}{\lambda}<0$; then:
	\begin{enumerate}
	\item Whenever $P_-\in T_-^1$ there is a unique $P_+\in T^1_+$ 
		such that the pair $(P_-,P_+)$ satisfies the Rankine-Hugoniot relations \eq{rh_VC}.
	\item Conversely, whenever $P_+\in T^1_+$  there is a unique $P_-\in T^1_-$ 
		such that the pair $(P_-,P_+)$ satisfies the Rankine-Hugoniot relations \eq{rh_VC}.
	\end{enumerate}
	In either case, $(P_-,P_+)$  is an admissible 
	self-similar 1-shock with similarity parameter $\lambda$, left state $P_-$, 
	right state $P_+$, and propagating along $x=\bar\xi t^\frac{1}{\lambda}$.
\end{proposition}

\begin{figure}[ht]
	\centering
	\includegraphics[width=8cm,height=8cm]{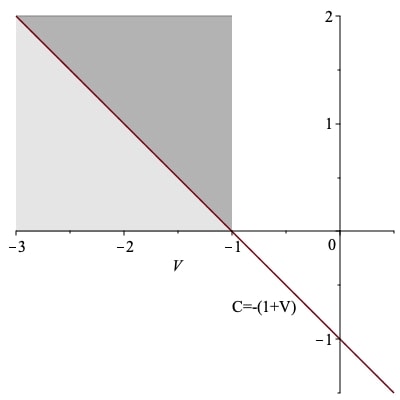}
	\caption{The regions $T^2_-=\{(V,C): 0<-(1+V)<C\}$ (dark grey) 
	and $T^2_+:=\{(V,C): 0<C<-(1+V)\}$ (light grey) of 
	possible left and right states, respectively, of an admissible 
	2-shock when $\frac{\bar \xi}{\lambda}<0$.}\label{T2_fig}
\end{figure}

Finally, for 2-shocks we define the regions
\beq\label{T2}
	T^2_-:=\{(V,C): 0<-(1+V)<C\},\qquad
	T^2_+:=\{(V,C): 0<C<-(1+V)\};
\eeq
see Figure \ref{T2_fig}. 
Since $\frac{\bar \xi}{\lambda}<0$,
the points in $T^2_-$ ($T^2_+$, respectively) satisfy the leftmost (rightmost, respectively)
inequality in the entropy condition \eq{2_shock} for a 2-shock. Again, the proof of the following proposition
is similar to that of Proposition \ref{SS_1_shocks_with_xi/lambda>0}.
\begin{proposition}[Admissible self-similar 2-shocks with $\frac{\bar \xi}{\lambda}<0$]
\label{SS_2_shocks_with_xi/lambda<0}
	Assume $\frac{\bar \xi}{\lambda}<0$; then:
	\begin{enumerate}
	\item Whenever $P_-\in T_-^2$ there is a unique $P_+\in T^2_+$ 
		such that the pair $(P_-,P_+)$ satisfies the Rankine-Hugoniot relations \eq{rh_VC}.
	\item Conversely, whenever $P_+\in T^2_+$  there is a unique $P_-\in T^2_-$ 
		such that the pair $(P_-,P_+)$ satisfies the Rankine-Hugoniot relations \eq{rh_VC}.
	\end{enumerate}
	In either case, $(P_-,P_+)$  is an admissible 
	self-similar 2-shock with similarity parameter $\lambda$, left state $P_-$, 
	right state $P_+$, and propagating along $x=\bar\xi t^\frac{1}{\lambda}$.
\end{proposition}
\begin{remark}\label{unique_jump1}
	It is clear that the open sets $S_\pm^i$, 
	$T_\pm^i$, $i=1,2$, are disjoint and contains all points in the $(V,C)$-plane
	off the critical lines $L_\pm$. Thus, for each state 
	$\bar P\in\RR^2\smallsetminus(L_+\cup L_-)$
	there is a unique state, distinct from $\bar P$ and also located 
	off $L_+\cup L_-$, to which it can be joined via an 
	admissible self-similar shock. Whether the resulting shock is a 1- or a 2-shock, 
	and whether $\bar P$ is the left or the right state, depend on 
	which of the regions $S_\pm^i$, $T_\pm^i$, $i=1,2$, it belongs to.
\end{remark}
\begin{remark}\label{unique_jump2}
	The analysis above shows in particular that an admissible self-similar shock
	must cross one of the critical lines $L_\pm$ with $\sgn(C)$ unchanged.
\end{remark}

\subsection{Hugoniot locus of a trajectory}\label{hug_locus}
\begin{definition}\label{hug_loc_def}
	Given any part $\Gamma$ of a trajectory of \eq{CV_ode} parametrized by
	$\xi\mapsto (V(\xi),C(\xi))$. Assume that $\Gamma$ is located within one of 
	the open regions $S_\pm^i$, $T_\pm^i$, $i=1,2$, defined in Section \ref{rh_e_conds}.
	According to Remark \ref{unique_jump1}, for 
	each state $(V(\xi),C(\xi))$ there is a unique distinct state
	to which it can connect 
	via an admissible self-similar shock. As $\xi$ varies, these states trace 
	out a continuous curve which we refer to as the {\em Hugoniot locus of 
	$\Gamma$} and denote by $\mathrm{Hug}(\Gamma)$. 
\end{definition}
The following lemmata will be used to argue that some
of the data \eq{vac_ss_data} necessarily generate a shock wave 
emanating from the initial vacuum interface.
Lemma \ref{hug_1} addresses the case when a trajectory approaches 
$P_1$ while Lemma \ref{hug_2} concerns the behavior of $\mathrm{Hug}(\Gamma)$
when $\Gamma$ tends to infinity. 
\begin{lemma}\label{hug_1}
	If a trajectory $\Gamma$ of \eq{CV_ode} approaches one of 
	the critical lines $L_\pm$ at a point distinct from $P_1=(-1,0)$, 
	then $\mathrm{Hug}(\Gamma)$ approaches the same point. The same 
	holds if $\Gamma$ approaches $P_1$ with $C^2/(1+V)$ bounded. 
\end{lemma}
\begin{proof}
Assume for concreteness that the running point along the trajectory 
$\Gamma$ is the left state $(V_-,C_-)$, so that $(V_+,C_+)$ runs 
along $\mathrm{Hug}(\Gamma)$. By solving the first of the Rankine-Hugoniot 
relations in \eq{rh_VC} for $1+V_+$ and substituting into the second, 
we obtain the equation
\beq\label{hug_limit}
	\textstyle\frac{(1+V_-)^2}{C_-^2}z^{2(\ell+1)}
	-(\frac{(1+V_-)^2}{C_-^2}+\frac{1}{\gamma})z^{\ell+2}+\frac{1}{\gamma}=0
	\qquad\text{for $z:=|\frac{C_-}{C_+}|$.}
\eeq
If now $\Gamma$ tends to a point on $L_+\cup L_-$ different from $P_1$,
then $\frac{(1+V_-)^2}{C_-^2}=1$ there, and the corresponding limiting $z$-value 
therefore satisfies
\[f(z):=\textstyle z^{2(\ell+1)}
-(1+\frac{1}{\gamma})z^{\ell+2}+\frac{1}{\gamma}=0.\]
It is immediate to verify that $f(0)>0$, $f(1)=0$,  
$f'(z)=2(\ell+1)z^{\ell+1}(z^\ell-1)$, so that $z=1$
is the unique root of $f$, i.e., $|C_+|=|C_-|$. It follows from Remark 
\ref{unique_jump2} that $C_+=C_-$, and the first Rankine-Hugoniot 
relation \eq{rh_VC}${}_1$ then gives $V_+=V_-$. 

Next, assume $P_1=(-1,0)$ is approached with $C_-^2/(1+V_-)$ bounded as $V_-\to-1$.  
It follows from \eq{hug_limit} that the limiting $z$-value then satisfies $-z^{\ell+2}+1=0$,
so that again $z=1$. 
Since $C_-\to 0$ as $P_1$ is approached, this shows that also $C_+\to0$.
Finally, since $z=|\frac{C_-}{C_+}|\to 1$ and $V_-\to-1$ as $P_1$ is approached,
\eq{rh_VC}${}_1$ gives $V_+\to-1$ as well. Thus, $(V_+,C_+)\to P_1$ as well.
\end{proof}

\begin{lemma}\label{hug_2}
	Let $\Gamma$ be a trajectory of \eq{CV_ode} and assume
	$(V_+,C_+)\in\Gamma$ tends to infinity with limiting 
	slope $k$. Then the corresponding Hugoniot point $(V_-,C_-)\in\mathrm{Hug}(\Gamma)$ 
	tends to infinity with a limiting slope $\tilde k$.
	If $\Gamma$ tends to infinity in $T_+^2$, then
	$\mathrm{Hug}(\Gamma)$ tends to infinity in $T_-^2$ with limiting slope $\tilde k\geq-1$.
\end{lemma}
\begin{proof}
For later use we assume that the running point along the trajectory 
$\Gamma$ is the right state $(V_+,C_+)$ and $(V_-,C_-)$ runs along 
$\mathrm{Hug}(\Gamma)$. Proceeding as in the proof of Lemma \ref{hug_1}, the 
Rankine-Hugoniot conditions imply that 
\beq\label{hug_limit_2}
	\textstyle\frac{(1+V_+)^2}{C_+^2}w^{2(\ell+1)}
	-(\frac{(1+V_+)^2}{C_+^2}+\frac{1}{\gamma})w^{\ell+2}+\frac{1}{\gamma}=0
	\qquad\text{where $w:=\frac{|C_+|}{|C_-|}$.}
\eeq
If $\Gamma$ approaches infinity with a limiting slope $k$, i.e., $\frac{C_+}{V_+}\to k$
as $|V_+|, |C_+|\to \infty$, then \eq{hug_limit_2} gives that
the corresponding limiting $w$-value 
\[w_k=\lim_{|V_+|,|C_+|\to\infty}\textstyle\frac{|C_+|}{|C_-|},\] 
is a root of the function
\[g(w):=\textstyle \frac{1}{k^2} w^{2(\ell+1)}
-(\frac{1}{k^2}+\frac{1}{\gamma})w^{\ell+2}+\frac{1}{\gamma}.\]
It is immediate to verify that $g(0)>0$, $g(1)=0$, $g(w)\to\infty$ as $w\uparrow\infty$,
and that $g$ has a global minimum at $w=\bar w_k:=(\frac{\gamma+k^2}{\gamma+1})^{1/\ell}$.
Note that $\bar w_k\gtrless 1$ according to $k^2\gtrless1$. It follows that
$g$ has the unique root $w_k=1$ when $k^2=1$, and that 
$g$ has a unique root $w_k$ different from $1$ when $k^2\neq1$,
with $w_k\gtrless1$ for $k^2\gtrless1$. Thus
\[\textstyle\frac{C_+}{V_+}\to k\qquad\text{and}\qquad \textstyle\frac{C_+}{C_-}\to w_k\]
as $(V_+,C_+)$ tends to infinity along $\Gamma$.
It follows that $C_-\to\infty$ and from \eq{rh_VC}${}_1$ we deduce that
\[\textstyle\frac{C_-}{V_-}\to \tilde k:=kw_k^{-1-\ell},\]
showing that $(V_-,C_-)$ tends to infinity with limiting slope $\tilde k$.

Finally, consider the case when $\Gamma$ tends to infinity 
within the set $T_+^2$ with a limiting slope $k$. It follows from the 
definition of $T_+^2$ in \eq{T2} that $k\geq-1$. According to Proposition 
\ref{SS_2_shocks_with_xi/lambda<0} we have that $\mathrm{Hug}(\Gamma)$ 
is located within $T_-^2$ (see Figure \ref{T2_fig}), from which it follows that 
$\mathrm{Hug}(\Gamma)$ tends to infinity there with a slope $\tilde k\leq-1$. 
\end{proof}
\section{Resolution of Case (I): $0<\lambda<1$ and $1<\gamma<3$}
\label{I}
\begin{figure}[ht]
	\centering
	\includegraphics[width=8cm,height=8cm]{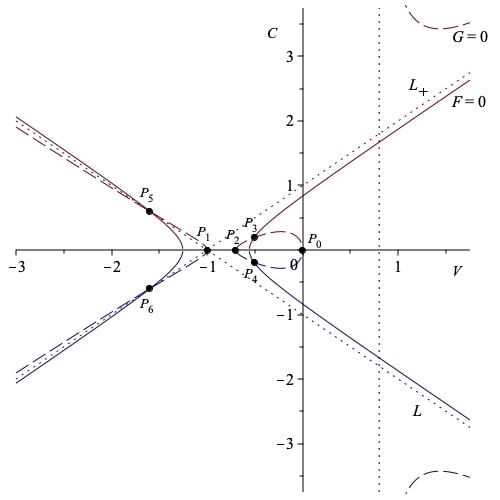}
	\caption{Case (I): The critical points $P_0$-$P_6$, the zero-level curves  
	of $F(V,C)$ (solid, including the $V$-axis) and $G(V,C)$ (dashed),
	the critical lines $L_\pm=\{C=\pm(1+V)\}$ (dotted), and the vertical asymptote 
	$V=V_*$ of $\{G=0\}$ (dotted). The parameters are $\lambda=0.7$ and
	$\gamma=1.75$.}\label{IA_fig_1a}
\end{figure}

\begin{figure}[ht]
	\centering
	\includegraphics[width=10cm,height=8.5cm]{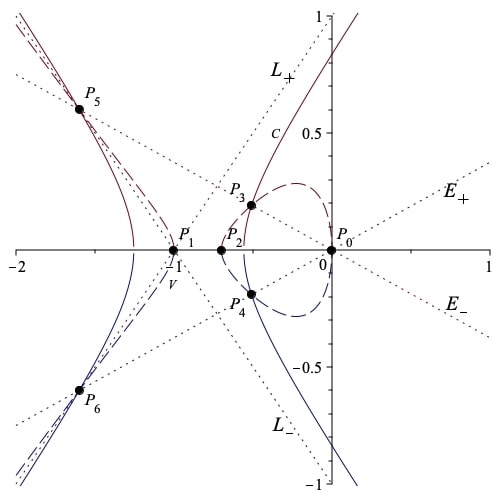}
	\caption{Case (I): Zoom-in of Figure \ref{IA_fig_1a} including 
	straight-line trajectories $E_\pm$.}\label{IA_fig_1b}
\end{figure}

We are now ready to describe the resolution of the initial value 
problems with vacuum data \eq{vac_ss_data}. This section addresses
the case in which the sound speed initially decays to zero in a H\"older manner
at the initial vacuum interface.
In this case the critical points $P_3,P_4$ lie within the cone $\mathcal K=\{|C|\leq|1+V|\}$: 
the condition \eq{P_3_in_K} reduces to $\lambda<\frac{\gamma+1}{3-\gamma}$,
which is met in this case. We also note that $-\lambda<V_3<0$ and $V_5<-1$.
See Figures \ref{IA_fig_1a}-\ref{IA_fig_1b} for a representative case.


\begin{figure}[ht]
	\centering
	\includegraphics[width=8cm,height=8cm]{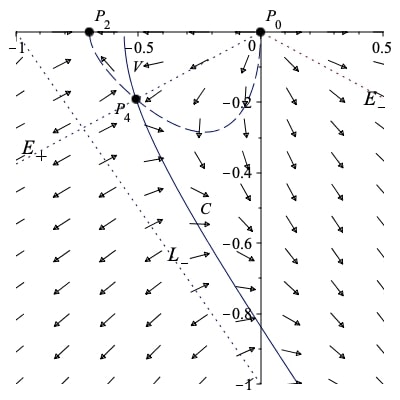}
	\caption{Case (I): Direction field plot of the similarity ODEs \eq{V_ode}-\eq{C_ode}
	in the lower half-plane near the origin.
	Arrows indicate flow direction as $\xi>0$ decreases. The critical line
	$L_-$ and the straight-line trajectories $E_\pm=\{C=\pm\ell^{-1}V\}$ 
	are dotted; the zero levels of $F$ and $G$ are solid and dashed, respectively.
	The parameters are the same as in Figures \ref{IA_fig_1a}-\ref{IA_fig_1b}.}
	\label{IA_fig_1c}
\end{figure}

\begin{figure}[ht]
	\centering
	\includegraphics[width=8cm,height=9cm]{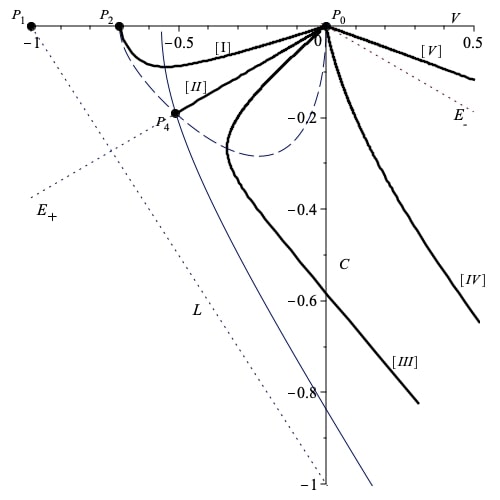}
	\caption{Case (I): $\Gamma_0$ trajectories $[I]$-$[V]$ displaying 
	different behavior depending on their slope at $P_0$.
	Trajectories $[I]$ and $[I\!I]$ connect $P_0$ directly to $P_2$ 
	and $P_4$, respectively (see Section \ref{IAa});
	the resulting vacuum interfaces remain stationary and no physical 
	singularity appears in these cases.
	The trajectories $[I\!I\!I]$-$[V]$ approach infinity in the fourth quadrant 
	and are continued by entering the second quadrant where they end 
	at $P_1$ (see Figures \ref{IA_fig_1e}-\ref{IA_fig_1g}); the resulting 
	vacuum interfaces propagate a physical singularity along a left-moving 
	interface. For $[I\!I\!I]$ and $[I\!V]$ the resulting Euler flow is continuous 
	(Section \ref{IAb}), while for $[V]$ it contains an admissible 2-shock 
	emanating from the initial vacuum interface (Section \ref{IAc}).  
	The parameters are as in Figures \ref{IA_fig_1a}-\ref{IA_fig_1c}.}
	\label{IA_fig_1d}
\end{figure}

According to the analysis in Section \ref{sim_odes} the $(V,C)$-trajectory $\Gamma_0$
selected by the initial data \eq{vac_ss_data} emanates from the 
origin $P_0$ with slope $\text{Ma}^{-1}=c_+/u_+$ as $\xi$ decreases from $\infty$. 
(Note that, since $c_+>0$, $\Gamma_0$ leaves the origin with nonzero slope.)
In the present case, with $\lambda>0$, \eq{approach_at_t=0} gives that the 
trajectory moves into the lower half-plane as $\xi$ decreases from $\infty$. Figure \ref{IA_fig_1c} 
displays the flow field of the original similarity ODEs \eq{V_ode}-\eq{C_ode} in the 
lower half-plane near the origin (with the same parameter 
values as in Figures \ref{IA_fig_1a}-\ref{IA_fig_1b}). Note that in Figure \ref{IA_fig_1c}  
the arrows indicate the direction of flow as $\xi>0$ {\em decreases}.

As is evident from Figure \ref{IA_fig_1c}, when the trajectory $\Gamma_0$ enters 
the third quadrant, its global behavior depends on whether it emanates from the 
origin above, along, or below the special straight-line trajectory 
$E_+=\{C=\ell^{-1}V\}$. Similarly, as we shall see below, when $\Gamma_0$ 
enters the fourth quadrant, its global behavior depends on its location relative to 
the other straight-line trajectory $E_-=\{C=-\ell^{-1}V\}$. The following sub-sections treat the various cases.

\subsection{Continuous flow for $\mathrm{Ma}\geq\ell$}\label{IAa}
First assume that $\text{Ma}>\ell$ so that the trajectory $\Gamma_0$ leaves 
the origin strictly above $E_+=\{C=\ell^{-1}V\}$.
Since $E_+$ is a trajectory, while $P_4$ is a saddle
(cf.\ Section \ref{P34}), $\Gamma_0$
then tends to the node $P_2$ as $\xi$ decreases. (Cf.\ Figure \ref{IA_fig_1c} and 
trajectory $[I]$ in Figure \ref{IA_fig_1d}.)
Its approach to $P_2$ dictates how the vacuum is attained from within 
the fluid in the resulting Euler flow.
This behavior is readily obtained by considering \eq{C_ode} as $(V,C)\to(-\lambda,0)$:
the leading order behavior is then given by
\beq\label{decay_to_vac_IA_1}
	\textstyle\frac{dC}{d\xi}=A\frac{C}{\xi},
	\qquad\text{where $A=\frac{\lambda\ell^{-1}}{1-\lambda}>0.$}
\eeq
Thus, as $P_2$ is approached, $C(\xi)$ decays to zero like $\xi^A$ as $\xi\downarrow0$.

The curves $t^{-\frac{1}{\lambda}}x=\xi\equiv constant>0$ foliate the 
entire quarter plane $\{x>0,t>0\}$ as $\xi$ ranges between $0$ and $\infty$. It 
follows that, in Case (I) with $\mathrm{Ma}>\ell$, the trajectory $\Gamma_0$ connects 
the origin $P_0$ to $P_2$ and defines a continuous Euler flow in all of $\{x>0,t>0\}$.
In particular, the vacuum interface is fixed at $x=0$, which corresponds to $\xi=0$. 
Since $C(\xi)\sim\xi^A$, with $A$ as in \eq{decay_to_vac_IA_1}, \eq{sim_var_c} gives
the following behavior for the sound speed as the vacuum is approached from within the 
fluid: 
\[c(t,x)\sim x^{1+A}\qquad \text{at any fixed time $t>0$, as $x\downarrow0$.}\]

As discussed in Section \ref{interpretation}, this behavior is reasonable 
on physical grounds: with $\text{Ma}=\frac{u_+}{c_+}>\ell$ the fluid is initially moving away
from the vacuum region sufficiently fast to counteract the positive pressure gradient, and
the vacuum interface remains at $x=0$ indefinitely. Furthermore, the rarefying effect of 
the initial motion is strong enough to immediately change the manner in which the vacuum is 
approached: initially the sound speed reaches vacuum in a H\"older manner 
($c(0,x)\sim x^{1-\lambda}$), while at positive times it does so in a $C^1$ manner 
($c(t,x)\sim x^{1+A}$, $A>0$).
 
%
%

For the limiting case that the $\Gamma_0$ emanates from the origin with slope 
equal to $\frac{1}{\text{Ma}}=\ell^{-1}$, the trajectory moves toward $P_4=(V_3,-C_3)$ 
along the straight-line trajectory $E_+$. (Cf.\ trajectory $[I\!I]$ in Figure \ref{IA_fig_1d}.) 
According to \eq{V_ode_along_E_+}
$V(\xi)$ near $P_4$ satisfies
\[\textstyle\frac{dV}{d\xi}\approx \frac{B}{\xi}(V-V_3),\]
where $B=-V_3/(V_3+\frac{\ell}{1+\ell})=\frac{\lambda}{1-\lambda}>0$. 
It follows that $P_4$ is approached with $\xi\downarrow0$. Thus, the resulting Euler
flow is defined on all of $\{x>0,t>0\}$ also in this particular case. However, differently
from the case with $\mathrm{Ma}>\ell$, we now have that $C(\xi)$ tends to a  
negative constant ($C_4=-C_3$) as $\xi\downarrow0$. According to \eq{sim_var_c} we therefore get that the 
sound speed decays {\em linearly} to zero as the the vacuum interface $\{x=0\}$ is 
approached at positive times. Again, there is an abrupt change in behavior along the 
interface from $t=0$ to $t>0$. 

This establishes part (a) of Theorem \ref{thm_1}.

\subsection{Continuous flow for $-\ell\leq\mathrm{Ma}<\ell$}\label{IAb}
For this range of $\text{Ma}=u_+/c_+$, as $\xi$ decreases from $\infty$, the trajectory 
$\Gamma_0$ emanates from the origin 
with either a positive slope $c_+/u_+>\ell^{-1}$ and moves into the third quadrant
of the $(V,C)$-plane and below $E_+$, or with a negative slope $c_+/u_+\in[-\infty,-\ell^{-1}]$ and
moves into the fourth quadrant below (or on) $E_-$. See Figure \ref{IA_fig_1d}.
In the former case, $\Gamma_0$ continues by crossing vertically the lower part of the loop 
of $\{G=0\}$, and then continues to the right 
and into the fourth quadrant. In Figure \ref{IA_fig_1d} these behaviors are displayed by 
trajectories $[I\!I\!I]$ and $[I\!V]$.

It follows that in all cases now under consideration, the corresponding solution 
of the similarity ODEs \eq{V_ode}-\eq{C_ode} is eventually located in the fourth 
quadrant, below (or, in the limiting case $\mathrm{Ma}=-\ell$, on) the straight-line trajectory 
$E_-=\{C=-\ell^{-1}V\}$, and above the branch of $\{F=0\}$ located there.
This implies that, as $\xi$ decreases further, the solution moves off to infinity
in the fourth quadrant; see Figures \ref{IA_fig_1c}-\ref{IA_fig_1d}.

We now make the following claims:
\begin{enumerate}
	\item $\Gamma_0$ tends to infinity in the fourth quadrant as $\xi\downarrow 0$, 
	asymptotically with a constant slope $k\in (-1,-\ell^{-1}]$. 
	\item The solution can be uniquely continued to negative 
	$\xi$-values by having it move in from infinity in the 
	second quadrant of the $(V,C)$-plane, with the same asymptotic 
	slope $k$. In the second quadrant the solution is located above (or on) 
	the straight-line trajectory 
	$E_-=\{C=-\ell^{-1}V\}$ and below the branch of $\{G=0\}$ located there,
	cf.\ Figures \ref{IA_fig_1a}-\ref{IA_fig_1b}.
	\item As $\xi<0$ decreases further, the solution approaches the 
	node at $P_5$, reaching it with a finite, negative $\xi$-value. The
	approach to $P_5$ is along the primary direction at $P_5$ (except in the 
	limiting case where the solution lies along $E_-$, which is the 
	secondary direction at $P_5$).
	\item There is a unique trajectory joining the node $P_5$ to the saddle 
	$P_1$, along which $P_1$ is reached with a finite, negative $\xi$-value $\xi_v$. 
\end{enumerate}

With these claims established (below), we will have obtained a solution $(V(\xi),C(\xi)$) of 
\eq{V_ode}-\eq{C_ode} defined for all $\xi\in(\xi_v,\infty)$. Finally, 
\eq{sim_var_u}-\eq{sim_var_c} yields the corresponding Euler flow defined by the 
data \eq{vac_ss_data}. The approach to $P_1$ by the ODE-solution $(V(\xi),C(\xi)$)
corresponds to approaching the vacuum in the resulting Euler flow, which therefore 
has a vacuum interface moving to the left along the curve $x=\xi_v t^{1/\lambda}$.
According to the analysis in Section \ref{P1}, a physical singularity is propagated along
the interface.

In order to argue for the claims above we change to coordinates 
\[\textstyle W:=\frac{1}{V}\qquad\text{and}\qquad Z:=\frac{1}{C},\]
so that approaching infinity in the $(V,C)$-plane corresponds to approaching 
the origin in the $(W,Z)$-plane.
In $(W,Z)$-variables the reduced similarity ODE \eq{CV_ode} takes the form
\beq\label{ZW_ode}
	\textstyle\frac{dZ}{dW}=\frac{Z(Z^2-W^2)+A(W,Z)}{W(Z^2-W^2)+B(W,Z)},
\eeq
where the higher order terms $A$ and $B$ are given by
\beq\label{d&e}
	A(W,Z)=WZ^3(2-k_1+\lambda W)\qquad\text{and}\qquad
	B(W,Z)=W^2(V_*W^2+Z^2(1+\lambda+\lambda W)).
\eeq
The {\em singular directions} of \eq{ZW_ode} at the origin are determined 
by the equation $Z^2-W^2=0$, or 
$\theta=\pm\frac{\pi}{4}, \pm\frac{3\pi}{4}$ (polar angle in the $(W,Z)$-plane).
According to Theorem 64 in \cite{algm} (pp.\ 331-332), there is a 
one-to-one correspondence between nonsingular directions $\theta$
and solutions of \eq{ZW_ode} which approaches the origin in the $(W,Z)$-plane
along the direction $\theta$. 

Translating back to $(V,C)$-coordinates we therefore have: for each 
slope $k\neq\pm1$ there is a unique pair of trajectories of the 
reduced similarity ODE \eq{CV_ode} which tends to infinity with $\frac{C}{V}\to k$
as $|V|,|C|\to\infty$. In particular, for $k<0$, $k\neq-1$ there are unique 
trajectories, one in the fourth quadrant and one in the second quadrant, tending 
to infinity with slope $k$. 

To see how $\xi$ behaves as a solution $(V(\xi),C(\xi))$ approach infinity, 
we substitute $C\sim kV$ into \eq{V_ode}-\eq{C_ode} to get 
\[\textstyle\frac{1}{V}\frac{dV}{d\xi}\sim-\frac{1}{\xi} 
\qquad\text{and}\qquad 
\frac{1}{C}\frac{dC}{d\xi}\sim-\frac{1}{\xi} 
\qquad\text{as $|V|,|C|\to\infty$.}\]
It follows that any solution of \eq{V_ode}-\eq{C_ode} which tends to infinity 
with slope $k\neq\pm1$ does so with $\xi$ approaching zero (either $\xi\to0+$ 
or $\xi\to0-$), and that 
\beq\label{app_to_infty}
	V(\xi)\sim\textstyle\frac{K_V^\pm}{\xi}, \quad 
	C(\xi)\sim\textstyle\frac{K_C^\pm}{\xi}\qquad\text{as $\xi\to0\pm$,}
\eeq
where $K_V^\pm,K_C^\pm$ are constants satisfying 
\beq\label{Ks}
	\textstyle\frac{K_C^\pm}{K_V^\pm}=k.
\eeq
For later reference we note that \eq{sim_var_u_alt}-\eq{sim_var_c_alt} 
together with \eq{app_to_infty} give the following limits from the right 
and left along $x=0$ in the resulting Euler flow:
\begin{align}
	& u(t,x)\to -\textstyle\frac{1}{\lambda}t^{\frac{1}{\lambda}-1} K_V^\pm=:u_\pm(t)
	\qquad\text{as $x\to0\pm,$}\label{u_x=0}\\
	& c(t,x)\to -\textstyle\frac{1}{\lambda}t^{\frac{1}{\lambda}-1}K_C^\pm=:c_\pm(t)
	\qquad\text{as $x\to0\pm.$}\label{c_x=0}
\end{align}

We now return to the trajectory $\Gamma_0$ considered earlier. 
The analysis above shows that it tends to infinity in the fourth 
quadrant as $\xi\downarrow 0+$ and with a certain slope $k$ 
between $-1$ (the limiting slope of the $\{F=0\}$-branch in the 
fourth quadrant) and $-\ell^{-1}$ (the constant slope of the 
trajectory $E_-$). This verifies Claim (1) above. According to 
\eq{Ks} we have $k=K_C^+/K_V^+$.

We then continue the solution $(V(\xi),C(\xi))$ to {\em negative} $\xi$-values by 
selecting the trajectory of \eq{CV_ode} which tends to infinity with the {\em same} 
slope $k$ in the second quadrant of the $(V,C)$-plane. (This corresponds to 
continuing the corresponding $(W,Z)$-trajectory through the origin in a 
$C^1$ manner.) According to the analysis above this selects a unique trajectory, 
denoted $\Gamma_0'$, of \eq{CV_ode}. Note that, as we move from the fourth 
quadrant with $\xi>0$ along $\Gamma_0$  to the second quadrant 
with $\xi<0$ along $\Gamma_0'$, the constraint \eq{sign_cond} remains satisfied 
since both $\xi$ and $C(\xi)$ change signs.

\begin{figure}[ht]
	\centering
	\includegraphics[width=8cm,height=8cm]{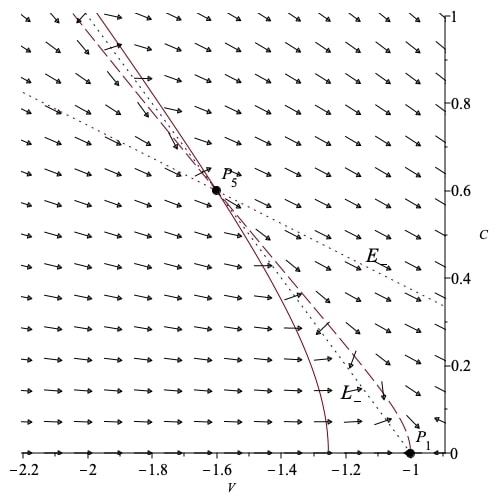}
	\caption{Case (I): Direction field plot of the similarity ODEs \eq{V_ode}-\eq{C_ode}
	in the upper half-plane near $P_5$.
	Arrows indicate flow direction as $\xi<0$ decreases. The critical line
	$L_-$ and the straight-line trajectory $E_-=\{C=-\ell^{-1}V\}$ 
	are dotted; the zero levels of $F$ and $G$ are solid and dashed, respectively.
	The parameters are the same as in Figures \ref{IA_fig_1a}-\ref{IA_fig_1b}.}
	\label{IA_fig_1e}
\end{figure}

\begin{figure}[ht]
	\centering
	\includegraphics[width=9cm,height=8.5cm]{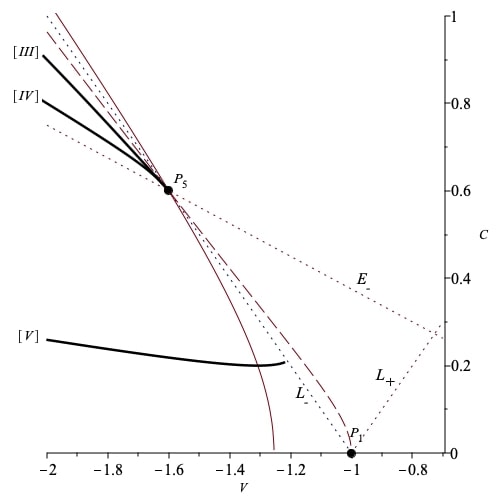}
	\caption{Case (I): $\Gamma_0'$ trajectories $[I\!I\!I]$-$[V]$ displaying different behavior depending on their slope 
	at infinity in the second quadrant.
	The parameters are as in Figures \ref{IA_fig_1a}-\ref{IA_fig_1b}.
	The critical lines $L_\pm$ and the straight-line trajectory $E_-$ are dotted. 
	The zero-levels $\{F=0\}$ and $\{G=0\}$ are solid and dashed, respectively.
	The trajectories $[I\!I\!I]$ and $[I\!V]$ come in from infinity between $E_-$ and $\{G=0\}$ 
	in the second quadrant and approach $P_5$. They will pass through $P_5$ and 
	end up at $P_1$ (see Figure \ref{IA_fig_1g}), defining corresponding Euler flows 
	that are continuous and with
	a physical singularity present along a left-moving interface (Section \ref{IAb}). 
	Finally, the trajectory $[V]$ comes in from infinity below $E_-$
	in the second quadrant; also this trajectory will end up at $P_1$ after jumping 
	across $L_-$ (see Figure \ref{IA_fig_1h}). The resulting flow contains a left-moving shock as
	well as a physical singularity located to the left of the shock (Section \ref{IAb}). } 
	\label{IA_fig_1f}
\end{figure}

\begin{remark}\label{diff_slopes}
	With the above construction we have that the solution $(V(\xi),C(\xi))$ 
	approaches infinity along $\Gamma_0'$ in 
         the second quadrant as $\xi\to0-$, and  with a limiting slope $K_C^-/K_V^-$.
         In particular, by insisting on the same slope at infinity for $\Gamma_0$ 
         and $\Gamma_0'$, we have $k=K_C^+/K_V^+=K_C^-/K_V^-=k'$.
         
	If instead $\Gamma_0'$ approached infinity with a slope $k'$ 
	{\em different} from $k$ it would follow from \eq{Ks} that either $K_V^+\neq K_V^-$ or
	$K_C^+\neq K_C^-$. In turn, \eq{u_x=0}-\eq{c_x=0} would imply the presence 
	of a stationary discontinuity along $x=0$ in the resulting Euler flow.
	While there is nothing wrong with this scenario per se, one can argue that
	for the particular self-similar solutions under consideration,  
	such a choice does not lead to a physically acceptable flow. 
	We omit the details of this argument which shows that $k'\neq k$ necessarily 
	leads to a solution which is either only partially defined or exhibits unphysical 
	behavior (infinite speed or propagation or entropy violating shocks).
\end{remark}

Since $\Gamma_0'$ tends to infinity with slope $k\in (-1,-\ell^{-1}]$ in 
the second quadrant, it follows that it must be located in the region $\mathcal R$
bounded below by the straight-line trajectory $E_-=\{C=-\ell^{-1}V\}$ 
and above by the $\{G=0\}$-branch in the second quadrant
(the latter having asymptotic slope $-1$ at infinity); see Figure \ref{IA_fig_1f}. This verifies 
Claim (2) above. 

Next, the phase portrait of \eq{sim_var_u}-\eq{sim_var_c} shows that the
region $\mathcal R$ is foliated by trajectories passing through the node $P_5$.
Thus, except for the limiting case when $\Gamma_0$ and $\Gamma_0'$ lie
along $E_-$, $\Gamma_0'$ approaches $P_5$ along its primary direction 
as $\xi$ decreases from zero. (Recall that the slope $-\ell^{-1}$ of $E_-$ 
gives the secondary direction at $P_5$, cf.\ \eq{prmr_scdr_slopes_5_2}.) 
The analysis of $P_5$ in Section \ref{P5} shows that it is reached as 
$\xi\downarrow \xi_5$, where $\xi_5$ finite and negative. This verifies Claim (3) above.

\begin{figure}[ht]
	\centering
	\includegraphics[width=12cm,height=8cm]{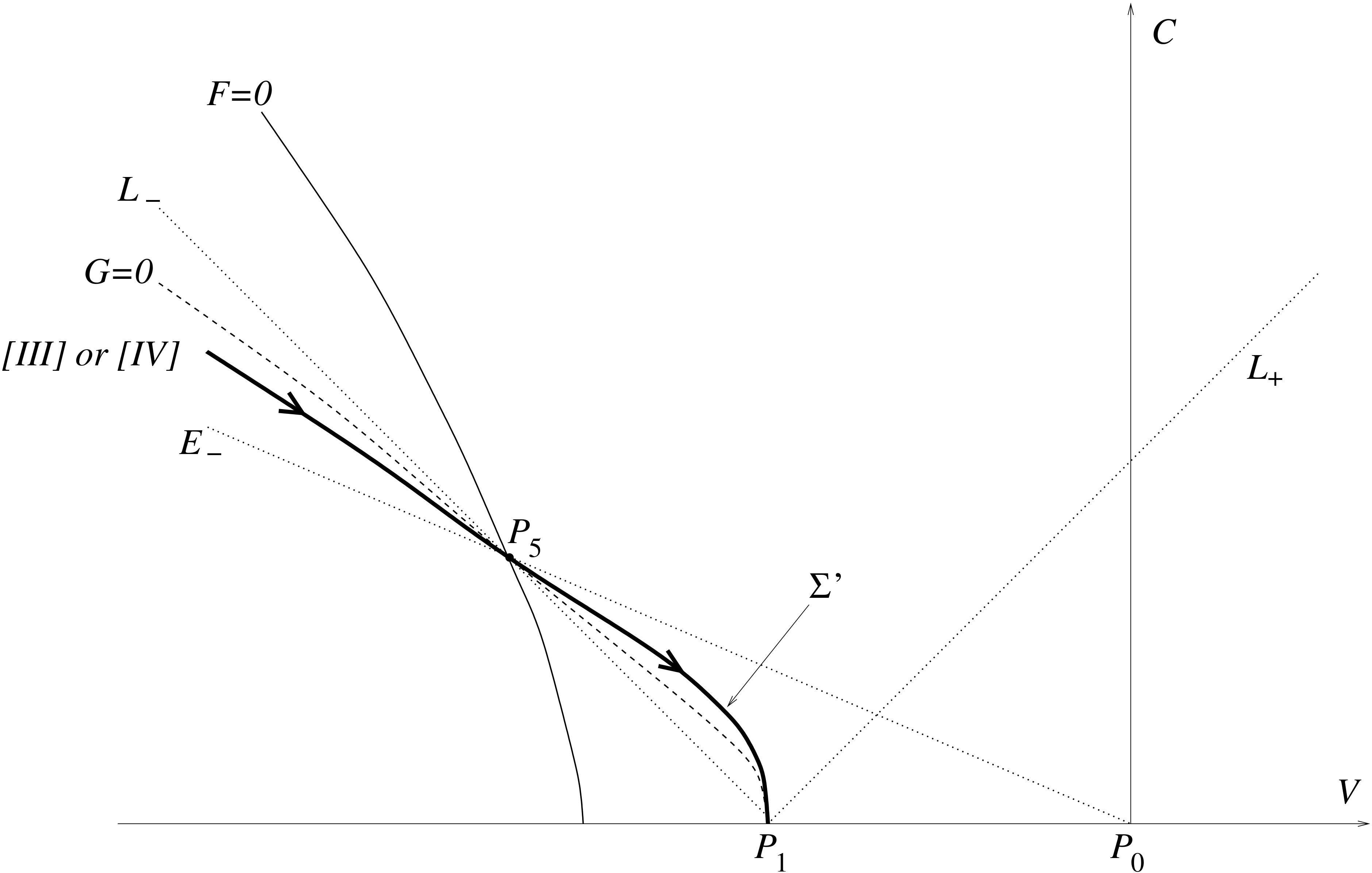}
	\caption{Case (I): Schematic figure of how trajectories $[I\!I\!I]$ and $[I\!V]$
	in Figure \ref{IA_fig_1f} pass through $P_5$ and continue to the saddle
	point $P_1$ along the separatrix $\Sigma'$. Arrows indicate 
	direction of motion as $\xi<0$ decreases.}
	\label{IA_fig_1g}
\end{figure}

To continue the solution as $\xi$ decreases beyond $\xi_5$,
we must select one of the infinitely many trajectories through 
the node $P_5$. An inspection of the phase portrait
shows that there is a unique solution ending up at the saddle $P_1$;
this is the separatrix denoted by $\Sigma'$ in Section \ref{P1}. 
Figure \ref{IA_fig_1g} provides a schematic picture of the situation.
Note that all other solutions of the similarity ODEs which leave 
$P_5$ with decreasing $\xi$ must necessarily run into one of the 
critical lines $L_\pm$, rendering them irrelevant for building a global
Euler flow. 
\begin{remark}
	Strictly speaking, the last statement requires elaboration. 
	The solutions leaving $P_5$ with decreasing $\xi$ and running into $L_+$
	(among them the one leaving along the straight-line trajectory $E_-$), 
	enter the region $T^1_+$ before reaching $L_+$. For these solutions there
	remains the possibility of performing an admissible jump across $L_+$ 
	into the region $T^1_-$, followed by a flow connecting to a vacuum state. 
	However, an analysis of the phase portrait in case (I) with $1<\gamma<3$
	reveals that one of three things can occur: the flow within $T^1_-$ brings 
	the trajectory under consideration to $P_0$, back to $L_+$, or to $P_3$. 
	According to analysis in Section \ref{outline} $P_0$ cannot serve as a 
	vacuum state. Also, approach to $L_+$ from within $T^1_-$ would occur 
	with a finite value of $\xi$, resulting in only a partially defined Euler flow 
	(note that performing another jump, now from $T^1_-$, is not admissible).
	Finally, an analysis similar to that for $P_0$ in Section \ref{outline} 
	shows that approach to $P_3$ must occur with $\xi$ tending to $-\infty$ 
	(we omit the details). However, this would imply an unphysical infinite speed 
	of propagation: the sound speed $c$, and hence the density, would be non-vanishing 
	along all of $\RR_x$ at any time $t>0$. 
	This justifies employing the unique node-saddle trajectory $\Sigma'$ to continue 
	the solution through $P_5$.
\end{remark}
According to the analysis of 
in Section \ref{P1}, $\Sigma'$ reaches $P_1$ with a finite $\xi$-value
$\xi_v<0$ (see \eq{lob_P1}). Finally, \eq{approach_to_vac_P1_2} shows that 
the resulting Euler flow connects to vacuum via a physical singularity.
This verifies Claim (4) and establishes part (b) of Theorem \ref{thm_1}.

\subsection{Discontinuous flow for $-\infty<\mathrm{Ma}<-\ell$}
\label{IAc}
For this range of $\text{Ma}=u_+/c_+$ the trajectory 
$\Gamma_0$ emanates from the origin with a negative slope 
$c_+/u_+\in(-\ell^{-1},0)$ as $\xi$ decreases 
from $\infty$. $\Gamma_0$ is then located in the region
between the $V$-axis and the straight-line trajectory
$E_-$ in the fourth quadrant. As in the previous case, 
$\Gamma_0$ moves off to infinity as $\xi\downarrow0$ with 
a certain slope $k\in(-\ell^{-1},0)$. In Figure \ref{IA_fig_1d} this 
behavior is displayed by trajectory $[V]$.

Arguing as above we get 
that the solution can be uniquely extended to negative
$\xi$-values by continuing along the unique trajectory 
$\Gamma_0'$ which approaches infinity with slope $k$ 
in the second quadrant. The trajectory $\Gamma_0'$ is then
necessarily located between the $V$-axis and the straight-line trajectory
$E_-$ in the second quadrant. In particular, $\Gamma_0'$ is 
located below the critical line $L_-$, see Figure \ref{IA_fig_1f}.

However, differently from the earlier cases, $\Gamma_0'$ now
cannot connect continuously to $P_1$:
$P_1$ is a saddle point which is approached only by the separatrix 
$\Sigma'$ which joins $P_5$ to $P_1$, and $\Sigma'$ is located above 
the critical line $L_-$. As is evident from Figures \ref{IA_fig_1e}-\ref{IA_fig_1f}, as $\xi$ 
decreases the solution moving along $\Gamma_0'$ must necessarily run into
$L_-$ at a point between $P_5$ and $P_1$ (and this occurs 
at a finite $\xi$-value).
Instead, we obtain a physically relevant solution by having an admissible 
jump occur from $\Gamma_0'$ (before it reaches $L_-$) to the separatrix 
$\Sigma'$.

\begin{figure}[ht]
	\centering
	\includegraphics[width=12cm,height=8cm]{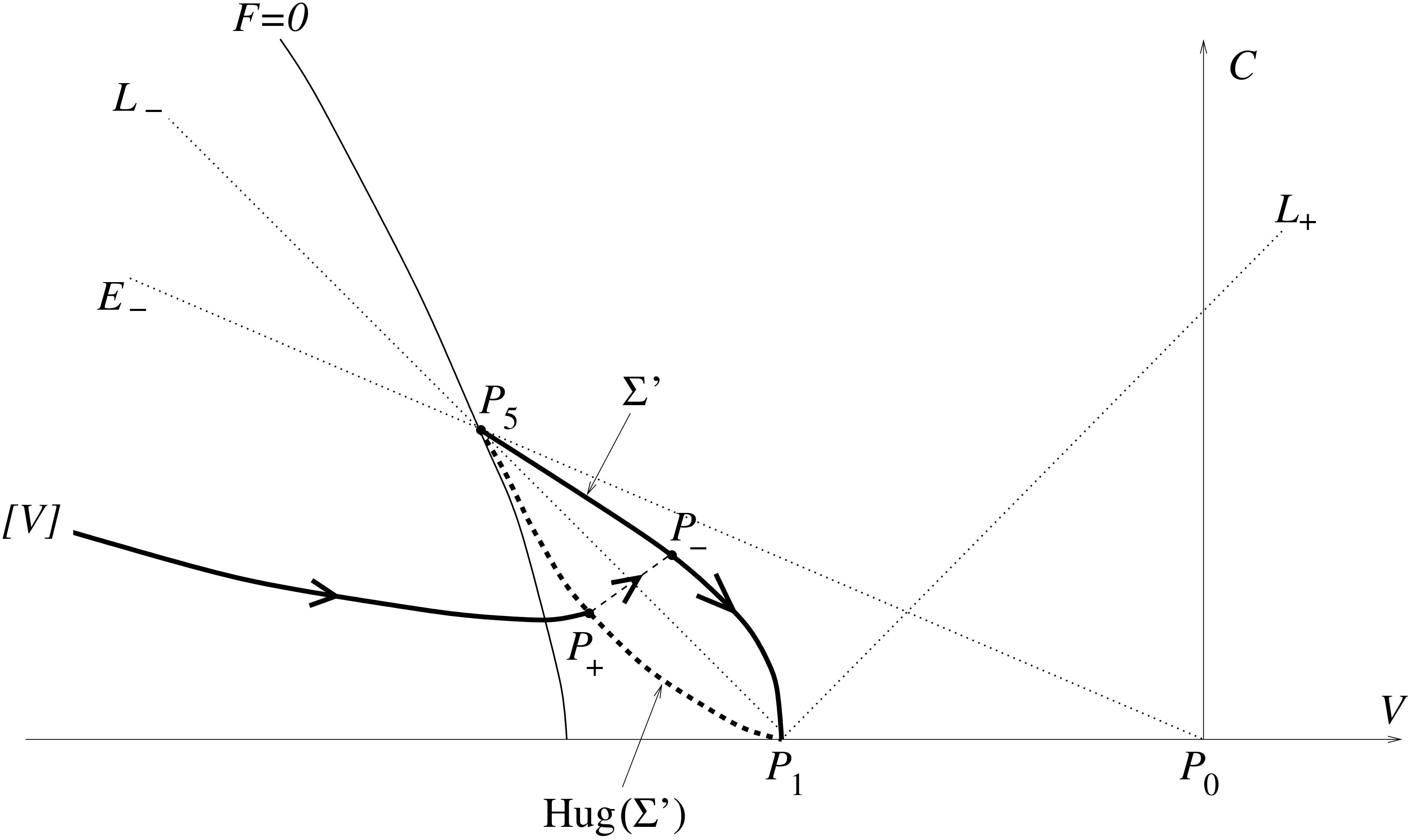}
	\caption{Case (I): Schematic figure of how trajectory $[V]$
	in Figure \ref{IA_fig_1f} jumps across $L_-$ from $P_+\in\mathrm{Hug}(\Sigma')$ to $P_-\in\Sigma'$, 
	and continues to the saddle point $P_1$ along the separatrix $\Sigma'$.
	Arrows indicate direction of motion as $\xi<0$ decreases.}
	\label{IA_fig_1h}
\end{figure}

We proceed to argue that this is always possible (see the schematic 
Figure \ref{IA_fig_1h}). Indeed, since 
both $P_1$ and $P_5$ are located on both $\Sigma'$ and $L_-$, 
we know from Lemma \ref{hug_1} that the Hugoniot locus
$\mathrm{Hug}(\Sigma')$ is a continuous curve which connects 
$P_1$ and $P_5$ (recall from Section \ref{P1} that $\Sigma'$ 
reaches $P_1$ with $C^2/(1+V)$ bounded, cf.\ \eq{P1_1st_separatrix}). 
Furthermore, since $\Sigma'$ is located within the region $T_-^2$,  it follows from part (1) 
of Proposition \ref{SS_2_shocks_with_xi/lambda<0} that its 
Hugoniot locus $\mathrm{Hug}(\Sigma')$ is located within the region 
$T_+^2$ (cf.\ Figure \ref{T2_fig}). Since $\Gamma_0'$ (trajectory $[V]$ in Figure \ref{IA_fig_1h})
approaches a point on $L_-$ between $P_1$ and $P_5$ from within $T_+^2$, 
it follows that $\Gamma_0'$ necessarily intersects $\mathrm{Hug}(\Sigma')$. 
Letting the point of intersection be denoted $P_+$, we get from 
part (2) of Proposition \ref{SS_2_shocks_with_xi/lambda<0} that 
$\Sigma'$ contains a corresponding point $P_-$ with the property 
that $(P_-,P_+)$  is an admissible self-similar 2-shock.

\begin{remark}\label{unique?}
	Numerical tests indicate that the intersection between 
	$\Gamma_0'$ and $\mathrm{Hug}(\Sigma')$ is unique; 
	however, we have not been able to prove this. If there 
	are multiple points of intersection, any one will work for
	our construction.
\end{remark}

Let $\xi_s$ denote the $\xi$-value for which the solution along $\Gamma_0'$
passes through $P_+$. Once the solution has jumped from $P_+\in\Gamma_0'$ 
to $P_-\in\Sigma'$ the analysis is as in the previous case: the solution moves 
on along $\Sigma'$ with $\xi$-values decreasing from $\xi_s$ and reaches 
$P_1$ with a finite $\xi$-value $\xi_v<\xi_s<0$. In particular, this gives the 
same type of behavior along the vacuum interface in the resulting Euler flow 
as in Section \ref{IAb}.

This establishes part (c) of Theorem \ref{thm_1}, and concludes the
proof of Theorem \ref{thm_1}.

\section{Resolution of Case (II): $\lambda< 0$ and $1<\gamma<3$}\label{II}
This section addresses the case where the similarity parameter $\lambda$ 
takes a negative value, i.e., when the initial sound speed $c_0(x)=c_+x^{1-\lambda}$ 
decays to zero in a $C^1$ manner at the initial vacuum interface $\{x=0\}$.
As in Section \ref{I} the adiabatic constant is restricted to $1<\gamma<3$.

\begin{figure}[ht]
	\centering
	\includegraphics[width=8cm,height=8cm]{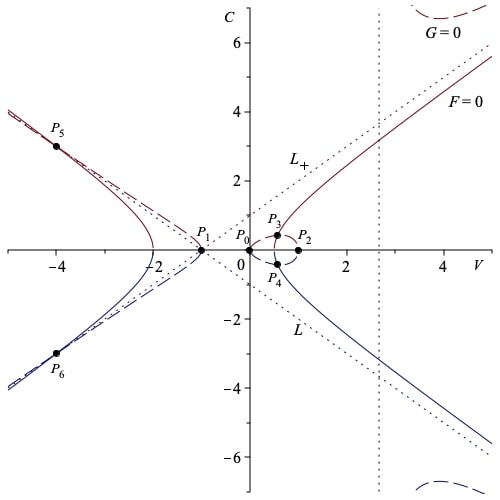}
	\caption{Case (II): The critical points $P_0$-$P_6$, the zero-level curves  
	of $F(V,C)$ (solid, including the $V$-axis) and $G(V,C)$ (dashed),
	the critical lines $L_\pm=\{C=\pm(1+V)\}$ (dotted), 
	and the vertical asymptote $V=V_*$ of $\{G=0\}$ (dotted). The parameters are 
	$\gamma=2.5$ and $\lambda=-1$.}\label{IIA_fig_1}
\end{figure}

The critical points $P_3$ and $P_4$ again belong to the 
cone $\mathcal K=\{|C|<|1+V|\}$ but are now located in the right half-plane. 
Figure \ref{IIA_fig_1} displays the critical points $P_0$-$P_6$, 
the zero levels of $F$ and $G$, together 
with the critical lines $L_\pm$ in a representative case.

According to the analysis in Section \ref{outline}, the trajectory $\Gamma_0$
selected by the initial data \eq{vac_ss_data} emanates from the 
origin with slope $1/\text{Ma}=c_+/u_+$ with $\xi$ 
increasing from $0$. Since $\lambda<0$, \eq{approach_at_t=0} 
gives that the trajectory moves into the upper half-plane as $\xi$ 
increases from zero. Figure \ref{IIA_fig_2} displays the flow field 
of the similarity ODEs \eq{V_ode}-\eq{C_ode} in the 
upper half-plane near the origin (the parameters are as in 
Figure \ref{IIA_fig_1}; note that the arrows are in the direction 
of increasing $\xi$). 
As in Case (I) (Section \ref{I}), the global behavior of $\Gamma_0$ 
depends on its initial slope at the origin.

\subsection{Continuous flow for $\ell\leq\mathrm{Ma}<\infty$}
\label{IIa}
\begin{figure}[ht]
	\centering
	\includegraphics[width=8cm,height=8cm]{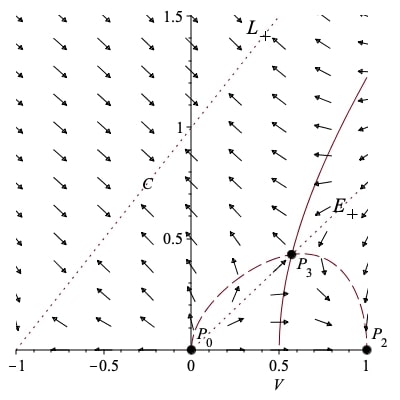}
	\caption{Case (II): Direction field plot of the similarity ODEs \eq{V_ode}-\eq{C_ode}
	in upper half-plane near the origin.
	Arrows indicate flow direction as $\xi>0$ increases. The critical line
	$L_+$ and the straight-line trajectory $E_+=\{C=\ell^{-1}V\}$ 
	are dotted; the zero levels of $F$ and $G$ are solid and dashed, respectively.
	The parameters are the same as in Figure \ref{IIA_fig_1}.}
	\label{IIA_fig_2}
\end{figure}

Consider first the case of strict inequality $\ell<\mathrm{Ma}<\infty$
so that $\Gamma_0$ leaves the origin with a slope $1/\text{Ma}$ strictly 
between $0$ and $\ell^{-1}$.
Since $E_+=\{C=\ell^{-1}V\}$ is a trajectory, while $P_3$ is a 
saddle (cf.\ Section \ref{P34}), any trajectory leaving the origin 
with a slope strictly between $0$ and $\ell^{-1}$ must approach the 
node $P_2$; see Figure \ref{IIA_fig_2}. 
Furthermore, by using the linearization of \eq{CV_ode} at $P_2$
(cf.\ \eq{P2_linzn}) together with the similarity ODEs \eq{V_ode}-\eq{C_ode}, 
one verifies that $P_2$ is necessarily reached as $\xi\uparrow\infty$. 
(The exact manner of approach depends on $\gamma\gtrless2$.) 

The curves 
$t^{-\frac{1}{\lambda}}x=\xi\equiv constant>0$ 
foliate the entire quarter plane $\{x>0,t>0\}$ as $\xi$ 
ranges from $0$ to $\infty$. It follows that each 
trajectory of \eq{V_ode}-\eq{C_ode} which 
leaves the origin with a slope between $0$ 
and $\ell^{-1}$ defines an Euler flow in all of 
$\{x>0,t>0\}$.

It remains to determine the behavior of the 
resulting flows when the boundary
$\{x=0, t>0\}$ is approached with $x\downarrow 0$. 
For $t>0$ fixed and $\lambda<0$, $x\downarrow 0$ 
corresponds to $\xi\downarrow 0$, and it follows 
from \eq{t=0_1} that to leading order 
\[C(\xi)\sim \xi^{-\lambda}\qquad\text{as $\xi\downarrow 0$.}\]
Therefore, by \eq{sim_var_c},
\[c(t,x)=\textstyle\frac{1}{\lambda}\frac{x}{t} C(\xi)\sim x^{1-\lambda}
\qquad\text{as $x\downarrow 0$ with $t>0$ fixed.}\]
This shows that the boundary 
curve $\{x=0, t>0\}$ is the vacuum interface in the resulting flow. Furthermore, at 
each time $t>0$, the sound speed decays to zero at the same rate $x^{1-\lambda}$ 
as it did initially. As discussed in Section \ref{interpretation}, this behavior is reasonable 
on physical grounds: with $\text{Ma}=u_+/c_+>\ell$ the fluid is initially 
moving away from the vacuum region sufficiently fast to counteract the positive 
pressure gradient, and the vacuum interface remains at $x=0$ indefinitely. 

For the limiting case that the trajectory leaves the origin with slope equal to $\ell^{-1}$,
the trajectory moves toward $P_3$ along the straight-line $E_+=\{C=\ell^{-1}V\}$.
According to \eq{V_ode_along_E_+} the behavior of $V(\xi)$ near $P_3$
satisfies
\[\textstyle\frac{dV}{d\xi}\approx \frac{B}{\xi}(V-V_3),\]
where $B=\frac{\lambda}{1-\lambda}<0$, and 
it follows that $P_3$ is approached with $\xi\uparrow\infty$. Thus, the resulting Euler
flow is defined on all of $\{x>0,t>0\}$ also in this particular case, and displays the 
same decay $c(t,x)\sim x^{1-\lambda}$ along the vacuum interface $\{x=0\}$ 
at all times $t\geq 0$.

\subsection{Partial flow for $\mathrm{Ma}<\ell$}
\label{IIb}
In this case the situation is different: the $\Gamma_0$ trajectory 
selected by the initial data \eq{vac_ss_data} leaves the origin with 
slope $\mathrm{Ma}^{-1}\in(\ell^{-1},\infty)\cup[-\infty,0)$. As is evident from Figure \ref{IIA_fig_2}, 
in these cases $\Gamma_0$ necessarily approaches the critical line 
$L_+=\{C=1+V\}$ at some point with $C>0$. It is readily verified that 
this occurs at a finite, positive $\xi$-value.
(It follows from \eq{non_obvious_reln} that all such trajectories approach 
$L_+$ in a north-west direction with slope $-\ell^{-1}$.)
In the present case with $\lambda< 0$ and $1<\gamma<3$, there is no triple 
point along this portion of $L_+$, ruling out the possibility of $\Gamma_0$ 
crossing $L_+$ in a continuous manner.

The only option for the trajectory to avoid running into $L_+$ would be to 
jump across it. As detailed in Section \ref{xi/lambda<0} (see Figure 
\ref{T1_fig}), an admissible 1-shock would do this. However, the trajectory 
would then necessarily jump to a state within $T^1_+:=\{(V,C): 0<1+V<C\}$
(cf.\ Proposition \ref{SS_1_shocks_with_xi/lambda<0}), and the trajectory 
would again flow toward $L_+$ (now approaching it in a south-east 
direction with slope $-\ell^{-1}$). 
Again this approach will occur at a finite $\xi$-value. 
Finally, it follows from the analysis in Section \ref{ss_grps_rh_e_conds} 
that no further admissible jump can be made from the set $T^1_+$. 

The upshot is that, in Case (II), whenever the data \eq{vac_ss_data} are such that 
the trajectory leaves the origin with a slope in $(\ell^{-1},\infty]\cup(-\infty,0)$, 
we obtain only a partially defined Euler flow within in a region of the form 
$\{(x,t)\,|\, t^{-\frac{1}{\lambda}}x<\xi_*\}$, where $0<\xi_*<\infty$.

As argued in Section \ref{interpretation}, this scenario is reasonable on 
physical grounds.
With $\text{Ma}=u_+/c_+<\ell$ the fluid is initially moving either away
relatively slowly from the vacuum region or toward it. Either way, the rapidly 
increasing pressure gradient in this case immediately generates 
an infinitely strong wave coming in from $x=+\infty$, leaving the 
flow undefined in its wake.

We summarize our findings for Case (II):
\begin{theorem}\label{thm_2}
	Assume $\lambda<0$ and $1<\gamma<3$.
	Consider the initial value problem for the 1-d isentropic Euler system 
	\eq{u}-\eq{c} with initial vacuum data \eq{vac_ss_data}.
	Let the signed Mach number of the data be $\text{Ma}=\frac{u_+}{c_+}$
	and let $\ell=\frac{2}{\gamma-1}$.
	Then: 
	\begin{enumerate}
		\item For $\text{Ma}\geq\ell$ there exists a globally defined, 
		self-similar, and shock-free solution. The solution has a stationary 
		vacuum interface along $x=0$ where the sound speed decays to 
		zero at the same super-linear rate $\sim x^{1-\lambda}$ at all times $t\geq0$.
		\item For $\text{Ma}<\ell$ the self-similar solution is defined only on a part 
		of the $(x,t)$-plane. 
	\end{enumerate}
\end{theorem}

\section{Other cases}\label{other_cases}
This section briefly summarizes the situation for the
cases with $\gamma\geq 3$. Since the analysis is
similar to that in Sections \ref{I} and \ref{II} we omit
most of the details.
\subsection{The case $\gamma=3$} 
In this case the initial value problem for \eq{mass}-\eq{mom}
with vacuum initial data of the form \eq{vac_ss_data} displays
the same qualitative features as in Cases (I) and (II). 
The differences in the phase plane are that the critical (triple) points 
$P_5,P_6$ are absent and the straight-line trajectories $E_\pm$ 
have slopes $\pm1$.

First, for $\lambda<0$ these differences play no role and the analysis 
is the same as in Sections \ref{IIa}-\ref{IIb}.
Next, for $0<\lambda<1$ we get that all $\Gamma_0$-trajectories 
leaving $P_0$ with slopes 
$\text{Ma}^{-1}\in(1,\infty)\cup[-\infty,-1]$ tend to infinity in 
the fourth quadrant with $\xi\downarrow0$ and asymptotic slope 
$-1$. These are all continued into the second quadrant along the 
{\em same} trajectory $\Sigma'$ (the separatrix in $\{C>0\}$ 
of the saddle point $P_1$), reaching $P_1$ with $\xi\downarrow\xi_v<0$
as in Section \ref{IAb}.
If instead $\Gamma_0$ leaves $P_0$ with a slope 
$\text{Ma}^{-1}\in(-1,0)$ it approaches infinity in the fourth 
quadrant with a certain slope in $(-1,0)$, and it is then continued as 
$\Gamma_0'$ with this latter slope into the second quadrant. 
It follows that $\Gamma_0'$ comes in from infinity within the region 
$T_+^2$ (cf.\ Figure \ref{T2_fig}). If continued it will reach the critical line
$L_-$ at a finite, negative $\xi$-value. However, by combining Lemmas 
\ref{hug_1} and \ref{hug_2} we have that the Hugoniot locus
$\mathrm{Hug}(\Sigma')$ goes through $P_1$ and tends to infinity 
with slope $-1$ within $T_+^2$. It must therefore intersect 
$\Gamma_0'$. It follows that a complete trajectory joining $P_0$ 
to $P_1$ can be built by jumping from $\Gamma_0'$ to $\Sigma'$,
as in Section \ref{IAc}. Finally, if $\Gamma_0$ leaves $P_0$ with a 
slope $\text{Ma}^{-1}\in(0,1]$ the analysis is the same as that in 
Section \ref{IAa}: in this case $\Gamma_0$ connects $P_0$ directly 
to $P_2$ or to $P_4$, and the resulting Euler flow is globally continuous.

The upshot is that Theorems \ref{thm_1} and \ref{thm_2} apply verbatim
also in the case $\gamma=3$ (with $\ell=1$).

%

\subsection{The case $\gamma>3$} 
Also in this case the solutions of \eq{mass}-\eq{mom}
with vacuum initial data \eq{vac_ss_data} display the same
qualitative features as in Cases (I) and (II) treated in Sections \ref{I}-\ref{II} 
(i.e., when $1<\gamma<3$). However, the phase portrait of the similarity ODEs 
\eq{V_ode}-\eq{C_ode} is now somewhat different.
Specifically, with $\gamma>3$, the critical points $P_5,P_6$ are 
located in the right half-plane. Also, the primary and secondary 
directions at these points are interchanged as  
$\lambda$ passes through the value $\hat\lambda\in(0,1)$ 
(cf.\ \eq{lambda_hat}).

For $\lambda<0$ these differences are irrelevant and
we obtain as before a globally defined
flow whenever the trajectory $\Gamma_0$ selected by the initial data 
\eq{vac_ss_data} leaves the origin with slope $\text{Ma}^{-1}\in(0,\ell^{-1}]$
and moves into the upper half-plane. 
$\Gamma_0$ connects $P_0$ to either $P_2$ (when $0<\text{Ma}^{-1}<\ell^{-1}$)
or $P_3$ (when $\text{Ma}^{-1}=\ell^{-1}$) and defines a continuous Euler flow. 
Again, no globally defined flow appears possible if $\text{Ma}^{-1}>\ell^{-1}$;
cf.\ Sections \ref{IIa}-\ref{IIb}.

\begin{figure}[ht]
	\centering
	\includegraphics[width=8cm,height=8cm]{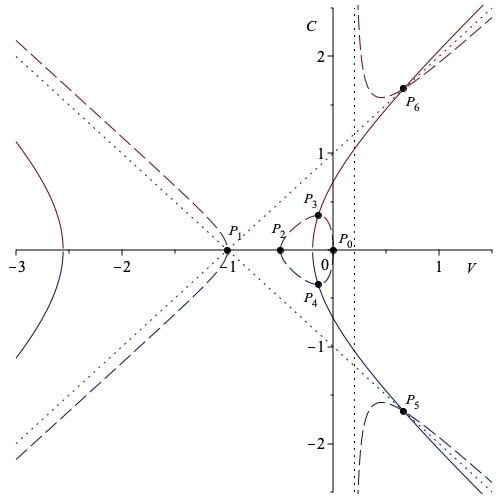}
	\caption{The case $\gamma>3$: The zero-level curves  
	of $F(V,C)$ (solid, including the $V$-axis) and $G(V,C)$ (dashed),
	together with the critical lines $L_\pm=\{C=\pm(1+V)\}$ (dotted) and the vertical asymptote 
	$V=V_*$ of $\{G=0\}$ (dotted). The parameters are $\gamma=6$ and
	$\lambda=0.5$.}\label{ib_fig}
\end{figure}

The situation for $0<\lambda<1$ is more involved. First,
the role played by the critical point $P_5$ 
in Case (I) is now played by $P_6=(V_5,-C_5)$,
which is located in the fourth quadrant.
The critical points $P_3,P_4$ again lie within the cone 
$\mathcal K=\{|C|\leq|1+V|\}$; see Figure \ref{ib_fig} 
for a representative case.
Also, Figure \ref{IA_fig_1c} for Case (I) still provides the 
correct qualitative features of the flow field of the similarity 
ODEs \eq{V_ode}-\eq{C_ode} in the lower half-plane near 
the origin. In particular, the trajectory $\Gamma_0$
selected by the initial data \eq{vac_ss_data} emanates from the 
origin $P_0$ with slope $1/\text{Ma}=c_+/u_+$, and 
moves into the lower half-plane  as $\xi$ decreases from $\infty$. 

It follows as in Section \ref{IAa} that when $\mathrm{Ma}\geq\ell$, 
the trajectory $\Gamma_0$ connects $P_0$ to either $P_2$
or $P_3$, and the result is a globally defined and continuous Euler flow.

When instead $\mathrm{Ma}<\ell$, and $\gamma>3$, a portion of the trajectories leaving $P_0$ 
will now pass through the triple point $P_6$ located in the fourth quadrant.
As demonstrated in Section \ref{P6}, the primary and secondary directions at $P_6$ 
are interchanged as the similarity parameter $\lambda$ crosses the value 
$\hat\lambda$. For certain ranges of the initial Mach number $\text{Ma}$ this
implies distinct types of behaviors in the corresponding Euler flow.
On the other hand, no fundamentally new feature appears in the 
solutions compared to Case (I): they will contain either a stationary or accelerating 
vacuum interface, and in the latter case a physical singularity is 
present, possibly together with an admissible 2-shock. 

Without going into the details we summarize our findings as follows. 
First assume $\hat\lambda<\lambda<1$. 
There is then a critical Mach number 
$\text{Ma}^*=\text{Ma}^*(\gamma,\lambda)<-\ell$
such that the following holds.

\begin{itemize}
	\item $-\ell\leq\mathrm{Ma}<\ell$: The resulting Euler flow is continuous and with a non-stationary 
	vacuum interface  along $x=t^\frac{1}{\lambda}\xi_v$ for a $\xi_v<0$. A physical singularity  
	is present along the vacuum interface at all times $t>0$. 
	\item $\mathrm{Ma^*}<\mathrm{Ma}<-\ell$: The resulting Euler flow contains a 2-shock 
	and a non-stationary vacuum interface. The vacuum interface moves to the left
	along $x=t^\frac{1}{\lambda}\xi_v$ for some $\xi_v<0$, while the 2-shock moves 
	to the right along $x=t^\frac{1}{\lambda}\xi_s$, where $\xi_s>0$. A physical singularity  
	is present along the vacuum interface at all times $t>0$. 
	\item $-\infty<\mathrm{Ma}\leq\mathrm{Ma^*}$: The resulting Euler flow contains a 
	2-shock and a non-stationary vacuum interface. The vacuum interface moves to the left
	along $x=t^\frac{1}{\lambda}\xi_v$ for some $\xi_v<0$, while the 2-shock moves 
	to the left along $x=t^\frac{1}{\lambda}\xi_s$, where $\xi_v<\xi_s<0$. A physical singularity  
	is present along the vacuum interface at all times $t>0$. 
\end{itemize}
Next consider the case that $0<\lambda<\hat\lambda$. 
There is then an additional critical Mach number 
$\text{Ma}^\circ=\text{Ma}^\circ(\gamma,\lambda)>\text{Ma}^*(\gamma,\lambda)$
such that the following holds.
\begin{itemize}
	\item $\text{Ma}^\circ\leq\mathrm{Ma}<\ell$: 
	The resulting Euler flow is continuous and with a non-stationary 
	vacuum interface  along $x=t^\frac{1}{\lambda}\xi_v$ for a $\xi_v<0$. A physical singularity  
	is present along the vacuum interface at all times $t>0$.  
	\item $\mathrm{Ma^*}<\mathrm{Ma}<\text{Ma}^\circ$: The resulting Euler flow contains a 2-shock 
	and a non-stationary vacuum interface. The vacuum interface moves to the left
	along $x=t^\frac{1}{\lambda}\xi_v$ for some $\xi_v<0$, while the 2-shock moves 
	to the right along $x=t^\frac{1}{\lambda}\xi_s$, where $\xi_s>0$. A physical singularity  
	is present along the vacuum interface at all times $t>0$. 
	\item $-\infty<\mathrm{Ma}\leq\mathrm{Ma^*}$: The resulting Euler flow contains a 
	2-shock and a non-stationary vacuum interface. The vacuum interface moves to the left
	along $x=t^\frac{1}{\lambda}\xi_v$ for some $\xi_v<0$, while the 2-shock moves 
	to the left along $x=t^\frac{1}{\lambda}\xi_s$, where $\xi_v<\xi_s<0$. A physical singularity  
	is present along the vacuum interface at all times $t>0$. 
\end{itemize}

\section*{Data availability statement}
The datasets used and/or analysed during the current study are available from the corresponding author on reasonable request.

\end{document}